\documentclass[12pt]{article}

\usepackage{amssymb,amsmath,amscd,amsthm}
\usepackage{graphicx,psfrag,epsfig}

\usepackage{graphicx}

\newtheorem{theorem}{Theorem}[section]

\newtheorem{lemma}[theorem]{Lemma}

\font\bbc=msbm10 scaled 1200
\newcommand{\E}{\mathbf{E}}
\newcommand{\R}{\mbox {\bbc R}}
\newcommand{\T}{\mbox {\bbc T}}
\newcommand{\Z}{\mbox {\bbc Z}}

\begin{document}

\title{Random Perturbations of 2-dimensional Hamiltonian Flows}

\author{ Leonid Koralov \\[1pt]
\normalsize Department of Mathematics\\[-4pt]
\normalsize Princeton University\\[-4pt]
\normalsize Princeton, NJ 08544\\[-4pt]
\normalsize koralov@math.princeton.edu\\[-4pt]
}

\date{}
\maketitle

\begin{abstract}
 We consider the motion  of a particle in a periodic two dimensional flow
perturbed by small (molecular) diffusion. The flow is generated by
a divergence free zero mean vector field. The long time behavior
corresponds to the behavior of the homogenized process - that is
diffusion process with the constant diffusion matrix (effective
diffusivity). We obtain the asymptotics of the effective
diffusivity when the molecular diffusion tends to zero.
\medskip
\end{abstract}

\section {Introduction}
\label{se1}

Consider the following stochastic differential equation
\begin{equation}
\label{e1} dX^\varepsilon_t = v(X^\varepsilon_t) dt + \sqrt{
\varepsilon} dW_t~,~~~ X^\varepsilon_t \in \R^2~.
\end{equation}
Here $v(x)$ is an incompressible
 periodic vector field, $W_t$ is a 2-dimensional
Brownian motion, and $\varepsilon $ (molecular diffusivity)
 is a small parameter. We further assume that the stream function $H(x_1,x_2)$,
such that
\[
\nabla^{\perp} H = (-H'_{x_2}, H'_{x_1}) = v~,
\]
is itself periodic in both variables, that is the integral of  $v$ over the
periodicity cell is zero.  For simplicity of notation  assume that the period of $H$ in each of the
variables is equal to one.

It is well known (see for example \cite{Fr}), that with $
\varepsilon$ fixed,  the solution of (\ref{e1}) scales like a
diffusion process with constant diffusion matrix when time goes to
infinity. More precisely, there exists the limit, called the
effective diffusivity,
\[
D^{ij}(\varepsilon) =  \lim_{t \rightarrow \infty}\E_{\lambda}
 \frac{X^{\varepsilon i}_t X^{\varepsilon j}_t}{t}~,~~~i,j=1,2~,
\]
where $i$ and $j$ are the coordinates and $\lambda$ is the initial
distribution of the process $X^{\varepsilon}_t$, which we can take
to be an arbitrary measure with compact support.
 The measure on $C([0,T], \R^2)$, induced by the process $\frac{1}{\sqrt{c}} X^\varepsilon_{ct}$,
converges weakly, when $c \rightarrow \infty$,
  to the measure induced by the diffusion process with constant matrix
$D(\varepsilon)$.

We are interested in the behavior of the effective diffusivity
when the molecular diffusion $\varepsilon$ tends to zero. Assume
that all
 the critical points of
$H$ are non degenerate. We distinguish two qualitatively different
cases,  depending on the structure of
 the stream lines of the flow given by $v(x)$.

 In the first case, there is a
 level set of $H$, which contains some of
 the saddle points, and  which  forms a lattice in $\R^2$, thus dividing
the plane into bounded sets, invariant under the flow. A standard
example of a cellular flow, which has been studied in several of
the papers cited below, is the flow with the stream function
$H(x_1, x_2) = \sin(2 \pi x_1) \sin(2 \pi x_2)$. In this
particular example the separatrices (the level sets of $H$
containing saddle points) form a rectangular lattice.

 In the second case, there is more than one unbounded level set of $H$ containing
 critical points, and thus
there are `open channels' in the flow, and some of the solutions
of the equation $x'(t) = v(x(t))$ go off to infinity. An example
of a flow with open channels is the flow with the stream function
$H(x_1, x_2)=\sin(2 \pi x_1) \sin(2 \pi x_2) + 10 \sin (2 \pi
x_2)$. Indeed, the horizontal axis $\{ x_2  = 0 \}$ is an
unbounded stream line of the flow.

\begin{figure}[htbp]
 \label{pic1}
  \begin{center}
    \begin{psfrags}
     \includegraphics[height=3in, width= 5.2in,angle=0]{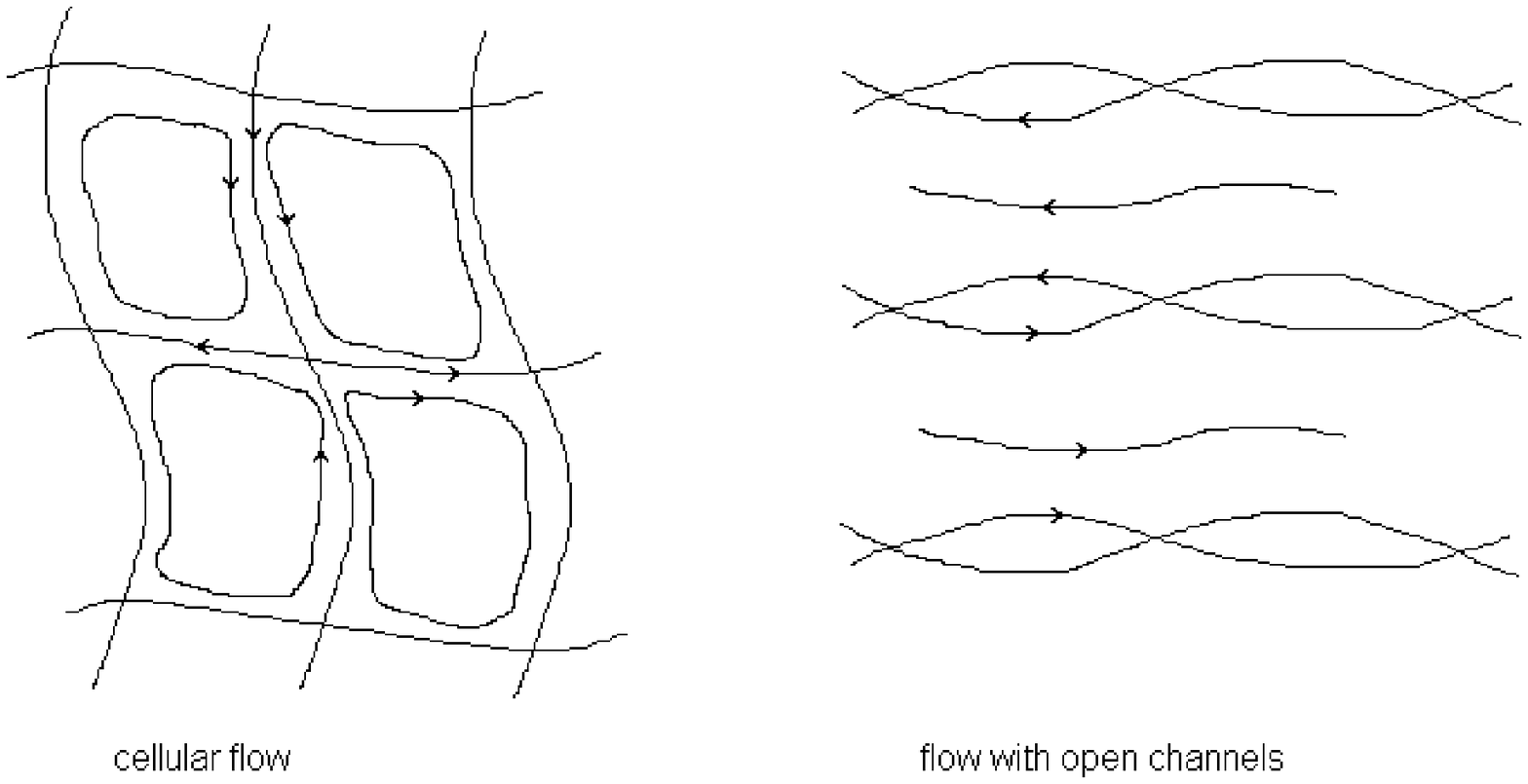}
    \end{psfrags}
  \end{center}
\end{figure}


 Since $v(x)$ is periodic,  we may consider $x'(t) = v(x(t))$ as the flow on the
torus. The torus is then  a union of the sepatatrices  and
   a finite number of open domains, bounded
by the separatrices, and invariant under the flow.

In \cite{FP} Fannjiang and Papanicolaou considered cellular flows
for which the separatrices form a rectangular lattice on $ \R^2$
and the stream function satisfies certain symmetry conditions.
They showed that in this case
\begin{equation}
\label{e5} D^{ij}(\varepsilon)= (d^{ij}+o(1)) {\sqrt{
\varepsilon}}~,~~~ {\rm as} ~~\varepsilon \rightarrow 0~,
\end{equation}
that is the effective diffusivity is enhanced by a factor of order
$\varepsilon^{-\frac{1}{2}}$ compared to case of the diffusion
process $\sqrt{\varepsilon} W_t$ without the advection term.
 Moreover, they found the constant matrix $d^{ij}$ explicitly.
Their proof is based on a variational principle applied to a
symmetric operator associated to the generator of the process
$X_t^{\varepsilon}$. In \cite{He} Heinze provided certain upper
and lower estimates on the effective diffusivity in the case of
cellular flows, for which the separatrices form a rectangular
lattice on $ \R^2$.

 There are earlier physical papers (\cite{Ch}, \cite{Shr}, \cite{So}),
 arguing that the asymptotics in (\ref{e5}) is true for particular flows.
 Our first result is the rigorous proof of this statement for general cellular
flows.

\begin{theorem}
\label{t1} Assume that an infinitely smooth periodic stream
function $H(x_1, x_2)$ defines a cellular flow, and that  its
critical points are nondegenerate.
 Then the asymptotics of the effective
diffusivity for the process (\ref{e1}) is given by (\ref{e5}).
\end{theorem}

Let ${\cal L}_p$ be the noncompact connected level set of $H$.
 This level set contains some
of the saddle points of $H$ and forms a lattice in $\R^2$.
Without loss of
 generality  we may assume that $H(x) = 0$ when $x \in {\cal L}_p$. The
corresponding level set on the torus will be denoted by ${\cal
L}$.

The process $X^\varepsilon_t$ consists of the `fast' part, which
is the periodic motion along the streamlines, and the `slow'
diffusion across them.
 The motion is almost periodic away
from the separatrices.
However,  once the trajectory is in a sufficiently small
neighborhood of the level set ${\cal L}_p$, it is likely to continue along it,
and may  go from cell to cell in a time much shorter than it would take the
`slow' diffusion to cover the same distance.

The rough outline of the proof of Theorem \ref{t1} is the
following. We introduce a Markov chain, which can be viewed as a
discrete time version of the  process $X^\varepsilon_t$. The state
space for the Markov chain
 is ${\cal L}$. Note, that due to the periodicity of $H$, the process $X^\varepsilon_t$
can be viewed as a process on the torus. In order to define the
transition probabilities, we introduce stopping times for the
process $X^\varepsilon_t$. The stopping time $\tau^\varepsilon_0$
is the first time when $X^\varepsilon_t$ hits ${\cal L}$, and
 $\tau^\varepsilon_n$ is defined as the first time after $\tau^\varepsilon_{n-1}$ when the
process $X^\varepsilon_t$ returns to ${\cal L}$, after having
traveled `past' a saddle
 point. The transition times of the Markov chain are random.

We show that the study of the asymptotics of the effective
diffusivity can be reduced to the study of the asymptotics of
transition probabilities and of the expectations of the transition
times for the Markov chain. The limit of the transition
probabilities as $\varepsilon \rightarrow 0$   is determined by
the behavior of the process $X^\varepsilon_t$
 in an arbitrarily small neighborhood of ${\cal L}$.
The asymptotics of the expectations of the transition times, on
the contrary, is determined by the event that  the trajectory of
$X^\varepsilon_t$ wanders away from the level
 set ${\cal L}$.

 In order to study the transition times we use
the results of Freidlin and Wentzell \cite{FW}.
 For a given stream function $H$ they  introduce a graph and
a mapping $g$  from the plane into the graph,
 such that each  connected level curve of
$H$ gets mapped into a point on the graph, with the level sets
 containing critical
points mapped into vertices. Then they demonstrate that the
process $g(X^\varepsilon_{t/\varepsilon})$ on the graph converges
to a limiting Markov process. The asymptotics of the expectations
of the transition times for our Markov chain is related to the
limiting process on the graph.

Now consider the flows with `open channels'. Assuming that the
channels are directed along the $x_1$ axis, we prove  that
\begin{equation}
\label{oc} D^{11}(\varepsilon)= (d^{11}+o(1)) \frac{1}{
\varepsilon}~, ~~~{\rm and}~~~ D^{22}(\varepsilon)= (d^{22}+o(1))
\varepsilon~,~~~ {\rm as} ~~\varepsilon \rightarrow 0~,
\end{equation}
that is the diffusion across the channels is not qualitatively
enhanced,
 compared with the process $ \sqrt{\varepsilon} W_t$.
The effective diffusivity  in the direction of the flow is
enhanced by a factor of order ${\varepsilon^{-2}}$.

\begin{theorem}
\label{t2} Assume that an infinitely smooth periodic stream
function $H(x_1, x_2)$ defines a flow with open channels, which
are directed along the $x_1$ axis, and that its critical points
are nondegenerate.
 Then the asymptotics of the effective
diffusivity for the process (\ref{e1}) is given by (\ref{oc}).
\end{theorem}

\noindent {\bf Remark} Since the matrix of effective diffusivity
is symmetric and positive definite, from Theorem \ref{t2} it
follows that the off-diagonal terms $D^{12}(\varepsilon) =
D^{21}(\varepsilon)$ are bounded uniformly in $\varepsilon$.
However we do not make a statement here on their asymptotic
behavior.
\\

 The proof of Theorem \ref{t2} is much simpler than that of
Theorem \ref{t1}, and is an easy application of the results of
Freidlin and Wentzell \cite{FW}.

The paper is organized as follows. In Section \ref{se2} we
describe the construction of the discrete time Markov chain
associated with a cellular flow, and relate the question of
effective diffusivity for the process $X^\varepsilon_t$ to the
study of transition probabilities and transition times for the
Markov chain. In Sections \ref{se3} and \ref{se4} we study of
transition probabilities and transition times respectively. In
Section \ref{se5} we prove Theorem \ref{t2}. In Section \ref{se6}
we prove several technical lemmas used in the previous sections.

\section{Construction of the Discrete Time Markov Chain}
\label{se2}

Consider the set $\{ x: |H(x)| < \varepsilon^{\alpha_1} \} $ on
the torus, with some $\frac{1}{4} < \alpha_1 < \frac{1}{2}$, and
let $V^\varepsilon$ be the connected component of this set, which
contains ${\cal L}$. Thus $V^\varepsilon$ is a thin tube around
${\cal L}$, whose width, however, is much larger than a typical
fluctuation of the process $H(X^\varepsilon_t)$ in fixed time (see
the picture below).

 Let $U_i,~  i = 1,...,n,$ be the connected components of $\T^2
\setminus {\cal L}$, and let $A_i,~ i=1,...,n,$ be the saddle
points, which belong to ${\cal L}$. While the numbers of the
connected components and of the saddles are the same for
topological reasons, their equality is not used in the proofs, and
the numbering of $U_i$'s is not related in any manner to that of
$A_i$'s. If there are points, which are carried to $A_j$ by the
flow $x'=v(x)$, and to $A_i$ by the flow $x'=-v(x)$, then the set
of such points (a subset of ${\cal L}$) is denoted by $\gamma(A_i,
A_j)$.
We assume (for the sake of simplicity of notation only) that
$\gamma(A_i, A_i)$ is empty, that is the separatrices do not form
"loops". In a neighborhood of each curve $\gamma(A_i, A_j)$ we may
consider a smooth change of coordinates $(x_1,x_2) \rightarrow (H,
\theta)$, where $\theta$ is defined by the conditions: $|\nabla
\theta| = |\nabla H|$ on $\gamma(A_i, A_j)$, and $\nabla \theta
\perp \nabla H$ (this way $\theta$ is defined up to multiplication
by $-1$ and up to an additive constant). The same change of
coordinates can be considered in $V^{\varepsilon}\cap
\overline{U}_k$. In this case $\theta \in [0,\int_{\partial U_k}
|\nabla H| dl ]$, with the end points of the interval identified.
Thus, if $A_i \in \partial U_k$, we define
\[
B(A_i, U_k) = \{ x \in V^\varepsilon \cap \overline{U}_k :
\theta(x) = \theta(A_i) \}.
\]
Let $B(A_i) = \bigcup_{k: A_i \in \partial U_k} B(A_i, U_k)$.

Define the stopping times $\sigma^\varepsilon_0 = 0, ~
\tau^\varepsilon_0 = \inf \{t: X^\varepsilon_t \in {\cal L} \} $.
Then $\sigma^\varepsilon_n, \tau^\varepsilon_n, ~ n \geq 1$ are
defined inductively as follows. Assume that
$X^\varepsilon_{\tau^\varepsilon_{n-1}}\in \gamma(A_i, A_j)$, and
$i \neq j$. Then
\[
\sigma^\varepsilon_{n} = \inf \{ t \geq \tau^\varepsilon_{n-1}:
X^\varepsilon_t \in \bigcup_{k \neq i} B(A_k) \bigcup \partial
V^\varepsilon \}.
\]
Thus, $\sigma^\varepsilon_n$ is the first time after
$\tau^\varepsilon_{n-1}$ that the process either exits
$V^\varepsilon$, or goes past a saddle point different from $A_i$.
Define $\tau^\varepsilon_n= \inf \{ t \geq \sigma^\varepsilon_{n}:
X^\varepsilon_t \in {\cal L} \}$. Let ${\cal L}^0 = {\cal L}
\setminus \{A_i, i=1,...,n\}$. Since almost every trajectory of
$X^\varepsilon_t$ does not contain any of the points $A_i$,
$X^\varepsilon_{\tau^\varepsilon_n}$ is a Markov chain with the
state space ${\cal L}^0$. The stopping times $\tau^\varepsilon_n$
are the consecutive times when the process $X^\varepsilon_t$ hits
the separatrix ${\cal L}$ after exiting $V^\varepsilon$ or after
having passed past a saddle point. Note, that the case, when a
point $x \in \gamma(A_i, A_j)$ travels, due to diffusion, against
the flow $v(x)$, and returns to $B(A_i)$, does not count as having
passed past a saddle point.

\begin{figure}[htbp]
 \label{pic2}
  \begin{center}
    \begin{psfrags}
     \includegraphics[height=3in, width= 5.2in,angle=0]{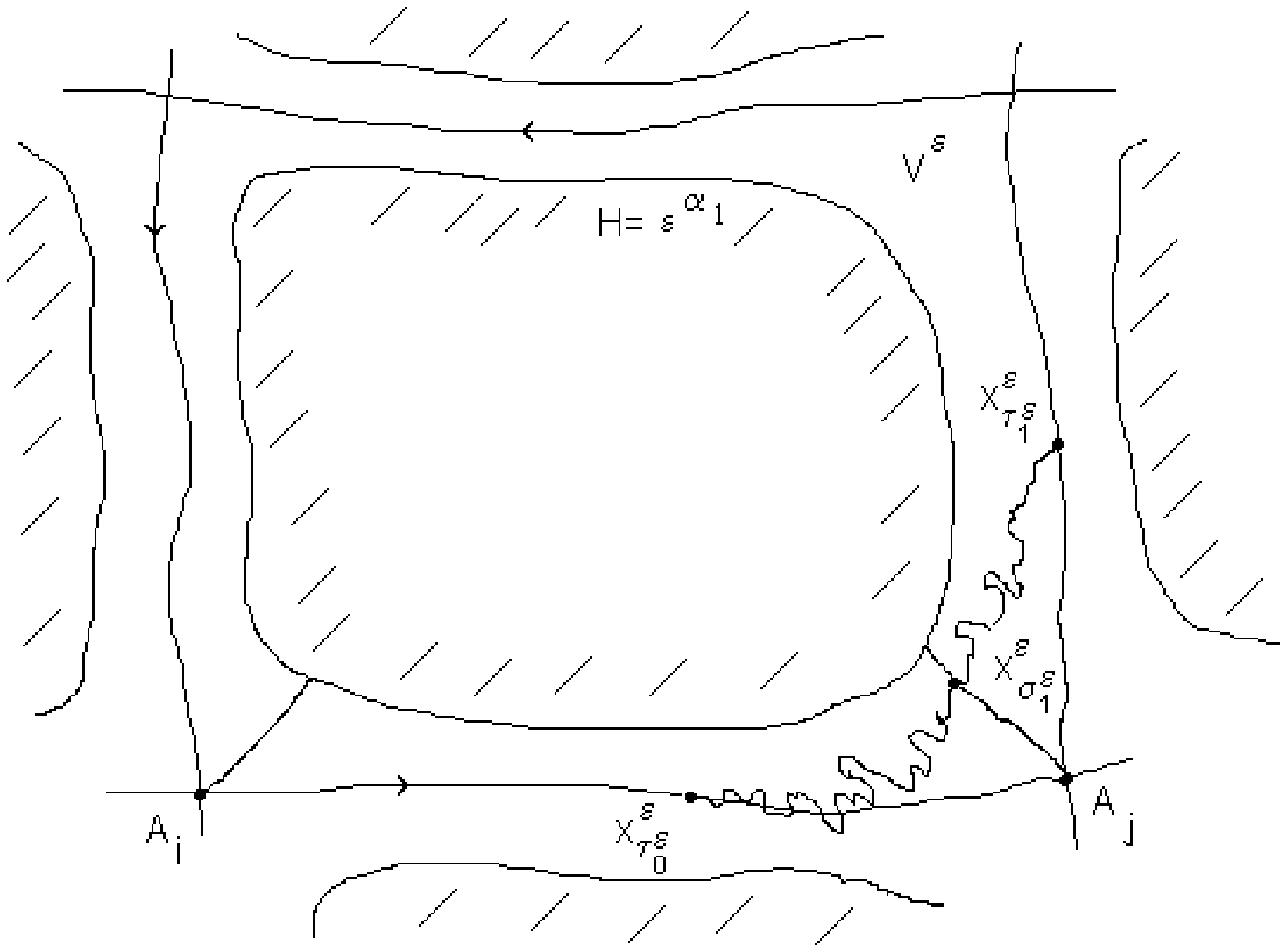}
    \end{psfrags}
  \end{center}
\end{figure}

It is not difficult to  see that
$X^\varepsilon_{\tau^\varepsilon_n}$ satisfies the Doeblin
condition, with the unique ergodic set. Therefore there exists a
unique invariant measure $\mu^{\varepsilon}(dy)$ on ${\cal L}^0$
(\cite{Doob}). Note that $(X^\varepsilon_{\tau^\varepsilon_n},
\tau^\varepsilon_n - \tau^\varepsilon_{n-1})$ also forms a Markov
chain with the state space ${\cal L}^0\times R_+ $, which
satisfies Doeblin condition for each $\varepsilon$, and has a
unique ergodic set . Let $\widetilde{p}^\varepsilon_x(dy, dt) $ be
the stochastic transition function  for this chain (it only
depends on the first component of the original point, as is
reflected in the notation). Then $ \widetilde{\mu}^\varepsilon
(dy, dt) = \int_{{\cal L}^0} \widetilde{p}^\varepsilon_x(dy, dt)
{\mu}^\varepsilon(dx) $ is the invariant measure. Since, due to
the presence of diffusion, the distributions of the transition
times have exponentially decreasing tails, we can apply the law of
large numbers to the function $(x,t) \rightarrow t$ on ${\cal
L}^0\times R_+ $ to obtain
\[
 \lim_{n \rightarrow \infty} \frac{\tau_n^{\varepsilon}}{n} =
\int\int_{{\cal L}_0\times R_+}t \widetilde{\mu}^\varepsilon(dy,
dt)= \]
\begin{equation}
\label{t}
 \int\int\int_{{\cal L}^0\times {\cal L}^0 \times R_+} t
\widetilde{p}^\varepsilon_x(dy, dt) \mu^\varepsilon(dx ) =
 \int_{{\cal L}^0} \E_x
\tau_1^{\varepsilon}  \mu^\varepsilon(dx)~~~~{\rm almost~surely}.
\end{equation}
In the arguments, which led to (\ref{t}), we considered
$X^{\varepsilon}_t$ as a process on the torus. In order to keep
track of the displacements of
$X^{\varepsilon}_{\tau^{\varepsilon}_n}$, as the process on
$\R^2$, we introduce another Markov chain, on the extended phase
space ${\cal L}^0\times \Z^2$. Let now $A_i, i=1,...$ be the
saddle points of $H$ on the plane. Then any $x_p \in \gamma(A_i,
A_j) \subset {\cal L}_p$ can be uniquely identified with a pair
$(x, z)$, where $x \in {\cal L}^0$ and $z = ([A^1_i], [A^2_i]) \in
\Z^2$.  ($[A^1_i]$ and $[A^2_i]$ are the integer parts of the
first and second coordinates of $A_i$). Thus, we have the mapping
$\phi: {\cal L}_p \setminus \{A_i, i=1,...\} \rightarrow {\cal
L}^0\times \Z^2$. Let $\phi_1$ and $\phi_2$ be the components of
this mapping. We define the Markov chain $Y^{\varepsilon}_n$ as
follows
\[
Y^{\varepsilon}_n =
(\phi_1(X^{\varepsilon}_{\tau^{\varepsilon}_n});
\phi_2(X^{\varepsilon}_{\tau^{\varepsilon}_n}) -
\phi_2(X^{\varepsilon}_{\tau^{\varepsilon}_{n-1}}))~.
\]
Note that the second component
$\phi_2(X^{\varepsilon}_{\tau^{\varepsilon}_n}) -
\phi_2(X^{\varepsilon}_{\tau^{\varepsilon}_{n-1}})$ almost surely
takes the values in some finite subset $S$  of $\Z^2$. Thus
$Y^{\varepsilon}_n$ is a Markov chain on ${\cal L}^0\times S$,
where $S=\{s_1,...,s_k \} \subset \Z^2$. It is not difficult to
see that $Y^{\varepsilon}_n$ satisfies Doeblin condition with the
unique ergodic set. Applying the law of large numbers to the
vector valued function $f(x,s) = s$, defined  on ${\cal L}^0\times
S$, we obtain that there exists $m \in \R^2$, such that
\[
\lim_{n \rightarrow \infty}
\frac{X^{\varepsilon}_{\tau^{\varepsilon}_n}}{n} = \lim_{n
\rightarrow \infty} \frac{\sum_{i=1}^n f(Y^{\varepsilon}_i)}{n}=
m~~~~{\rm almost~surely}~.
\]
Since $\frac{X^{\varepsilon}_t}{\sqrt{t}} \rightarrow N(0,
D(\varepsilon))$ in distribution, and due to (\ref{t}), we
conclude that $m=0$. Applying the central limit theorem to the
same function $f(x,s)$, we obtain that there exists a matrix
$d^{\varepsilon}$, such that
\[
\lim_{n \rightarrow \infty}
\frac{X^{\varepsilon}_{\tau^{\varepsilon}_n}}{\sqrt{n}} = \lim_{n
\rightarrow \infty} \frac{\sum_{i=1}^n
f(Y^{\varepsilon}_i)}{\sqrt{n} } = N(0, d^{\varepsilon})~~~~{\rm
in~distribution}~.
\]
Due to (\ref{t}), we have
\begin{equation}
\label{Dd} D(\varepsilon) = d^{\varepsilon} /\int_{{\cal L}^0}
\E_x \tau_1^{\varepsilon} d \mu^\varepsilon(x)~.
\end{equation}
In sections \ref{se3} and \ref{se4} we shall obtain the
asymptotics of the stochastic transition functions for the chain
$Y_n^\varepsilon$, and of the functions $\E_x
\tau_1^{\varepsilon}$. Assuming that this is accomplished, the
next lemma will allow us to obtain the asymptotics of the
effective diffusivity, using (\ref{Dd}).

First we introduce the notations and  formulate the assumptions
needed for the lemma. Let $M$ be a locally compact separable
metric space. Let $C_b(M)$ be  the set of bounded continuous
functions on $M$. Let $p_{\varepsilon}(x,dy), ~ 0 \leq \varepsilon
\leq \varepsilon^0$ be a family of stochastic transition functions
on $M$. Assume that

 (A) The family of measures $p_0(x,dy),~x \in K$ is tight for any
 compact set $K$.

 (B) $p_0(x,dy)$ is weakly Feller, that is $\int f(y) p_0(x,dy) \in C_b(M) $
if $f \in C_b(M)$.

 (C) For any $f \in C_b(M)$ and any compact $K \subset M$,
\[
\lim_{\varepsilon \rightarrow 0} \int f(y) p_{\varepsilon}(x,dy) =
\int f(y) p_0(x,dy)~~~~{\rm uniformly~in ~}x \in K~.
\]

 (D) There exist unique invariant measures
 $\mu^{\varepsilon}(dy)$. There exist $\lambda >0, c>0$, such that
 \[
 |p_{\varepsilon}^n (x, A) - \mu^{\varepsilon}(A)| \leq c e^{-
 \lambda n}~~~{\rm for ~ all~} x, A, \varepsilon~.
 \]
(That is  $p_{\varepsilon}$ are uniformly exponentially mixing).

Let $g \in C_b(M, \R^2)$ be such that $\int g
d\mu^{\varepsilon}=0$ for all $\varepsilon$. Let $Y_n^\varepsilon$
be the stationary Markov chain, with the stochastic transition
function $p_\varepsilon$. Since $Y^\varepsilon_i$ is exponentially
mixing, the central limit theorem can be applied to
$g(Y^\varepsilon_i)$, and thus $ \frac{\sum_{i=1}^n
g(Y^\varepsilon_i)}{\sqrt{n}}$ converges weakly as $n \rightarrow
\infty$ to a mean-zero Gaussian distribution. We denote the
covariance matrix of the limiting distribution by
 $d^\varepsilon(g)$.

\begin{lemma}
\label{l1} Suppose that assumptions (A)-(D) hold. Then \\ (a)
$\mu^{\varepsilon} \rightarrow \mu^0~~~{\rm weakly}.$\\ (b) If
$f^\varepsilon(x)$ are uniformly bounded, $f^0(x) \in C_b(M)$, and
$\lim_{\varepsilon \rightarrow 0} f^\varepsilon(x) = f^0(x)$
uniformly on any compact, then
\[
\lim_{\varepsilon \rightarrow 0} \int  f^\varepsilon d
\mu^\varepsilon = \int f^0 d\mu^0~.
\] \\
(c) If $g \in C_b(M, \R^2)$ is such that $\int g d\mu^\varepsilon
= 0$ for all $\varepsilon$, then
\[
d^\varepsilon(g) \rightarrow d^0(g)~.
\]
\end{lemma}
\noindent {\bf Proof:} From (A) it follows that for each $n$ the
family of measures $p_0^n(x,dy), x\in K$ is tight for any compact
set $K$. Let us assume that for a certain $n$, for any $f \in
C_b(M) $, uniformly on any compact set $K$  we have
\begin{equation}
\label{st} \int f(y) p_{\varepsilon}^n (x,dy) - \int f(y) p_{0}^n
(x,dy) \rightarrow 0~~{\rm as}~  \varepsilon \rightarrow 0~.
\end{equation}
Note that this  is true for $n =1$ by (C). Combining (\ref{st})
and the fact that $p_0^n(x,dy), x\in K$ is tight we obtain that
for any compact set $K$ and for any $\delta > 0$ there is a
compact set $K_1$ such that
\begin{equation}
\label{del}
p_{\varepsilon}^n(x, K_1)  > 1 - \delta, ~ x \in K
\end{equation}
for sufficiently small $\varepsilon$. Next we justify (\ref{st})
for $n+1$ instead of $n$.
\[
\int f(y) (p_{\varepsilon}^{n+1} (x,dy) -  p_{0}^{n+1} (x,dy))= \]
\[
\int f_1(y)(p_{\varepsilon}^{n} (x,dy) -  p_{0}^{n} (x,dy)) + \int
\int f(z) ( p_{\varepsilon}(y, dz) - p_0(y, dz) )
p_{\varepsilon}^n (x, dy) ~,
\]
where $f_1(x) = \int f(y) p_0(x, dy) \in C_b(M)$.  The first term
on the right hand side tends to zero by (\ref{st}), while the
second term tends to zero by (\ref{del}) as $\int f(z) (
p_{\varepsilon}(y, dz) - p_0(y, dz) )$ is bounded and tends to
zero uniformly on any compact. We therefore have established
(\ref{st}) for all $n$.

Let us prove part (a) of the lemma. Fix an arbitrary $x \in M$.
Then for $f \in C_b(M)$ we have
\[
\int  f (y) d\mu^{\varepsilon}(y) = \lim_{n \rightarrow \infty}
\int f(y) p_{\varepsilon}^n (x, dy)~,
\]
and the limit is uniform in $\varepsilon$ by (D). The weak
convergence of $\mu^{\varepsilon}$ to $\mu^0$ now follows from
(\ref{st}).

To prove part (b) we write
\[
\int f^{\varepsilon} d\mu^{\varepsilon} - \int f^0 d \mu^0 = (\int
f^0 d\mu^{\varepsilon} - \int f^0 d\mu^0) + (\int f^{\varepsilon}
d \mu^{\varepsilon} - \int f^0 d \mu^{\varepsilon} )~.
\]
The difference of the first two terms on the right hand side tends
to zero as $\mu^{\varepsilon} \rightarrow \mu^0$ weakly. The
difference of the last two terms  tends to zero since for any
$\delta
>0$ there is a  compact set $K$ for which  $ \mu^{\varepsilon}(K)
> 1 -\delta$ for sufficiently small $\varepsilon$ (since
$\mu^{\varepsilon} \rightarrow \mu^0$ weakly), and
$f^{\varepsilon} \rightarrow f^0$ on any compact set.

In order to prove part (c) of the lemma it is sufficient to
consider $g \in C_b(M, \R)$ (scalar valued). In this case
\[
d^{\varepsilon} (g) = \mathbf{E}[ (g(Y^{\varepsilon}_0))^2 +
2g(Y^{\varepsilon}_0)g(Y^{\varepsilon}_1) + 2
g(Y^{\varepsilon}_0)g(Y^{\varepsilon}_2)+...]~.
\]
Due to uniform mixing (D)
\[
| \mathbf{E}(g(Y^{\varepsilon}_0)g(Y^{\varepsilon}_n))| \leq
e^{-\gamma n}~,
\]
where $\gamma$ does not depend on $\varepsilon$.  In order to
prove that $d^{\varepsilon}(g) \rightarrow d^0(g)$ we therefore
only need to establish that
\begin{equation}
\label{nox}
 \mathbf{E}[g(Y^{\varepsilon}_0)g(Y^{\varepsilon}_n) -
g(Y^{0}_0)g(Y^{0}_n)] \rightarrow 0
\end{equation}
for any fixed $n$. The left hand side of (\ref{nox}) can be
written as
\[
\int g(x) g(y) p_{\varepsilon}^n (x, dy) d \mu^{\varepsilon} (x) -
\int g(x) g(y) p_{0}^n (x, dy) d \mu^{0} (x)~.
\]
Let $G_{\varepsilon}(x) = g(x) \int g(y) p_{\varepsilon}^n (x,
dy)$. Then $G_{\varepsilon}(x) \rightarrow G_0(x)$ uniformly on
any compact by (\ref{st}), and the conclusion follows by part (b).
\qed

\section{The Limit of the Transition Probabilities}
\label{se3}

In this section we shall identify the limit of the transition
probabilities for the chains $Y^{\varepsilon}_n$ on ${\cal
L}^0\times S$ and verify the conditions (A)-(D) for these chains.

Recall the $(H, \theta)$ coordinates which we may consider inside
each cell $U_k$ near its boundary, that is in $U_k \cap
V^{\varepsilon}$. Let
\[
h = \varepsilon^{-\frac{1}{2}}H~. \]
  In order to find the limit of the
transition probabilities we shall demonstrate that  in a small
neighborhood of ${\cal L}_p$ after a random change of time the
process $X_t^{\varepsilon}$ is well approximated by the process
$X_t$  with the generator $\frac{1}{2}\partial_{hh} +
\partial_{\theta}$ in $(h, \theta)$ coordinates.

Let $x \in \gamma(A_i, A_j)$ be a point on ${\cal L}_p$. We
introduce the stochastic transition function $p_0(x, dy), x,y \in
{\cal L}_p$ as follows:

Let $x, y \in \partial U_k$ (if $x$ and $y$ do not belong to the
boundary of the same cell, then $p_0(x, dy) = 0)$. Consider the
$(h, \theta)$ coordinates in $U_k \cap V^{\varepsilon}$, so that
$\theta(x) = 0$, and $\theta$ increases in the direction of the
flow. Since $\partial U_k$ is a closed contour, points with
coordinates $(h, \theta)$ and $(h, \theta + \int_{\partial U_k}
|\nabla H| dl) $ are identified. Let $\theta(A_j)$ and $\theta(y)$
belong to $(0, \int_{\partial U_k} |\nabla H| dl ]$, and consider
the process $X_t$ with the generator $\partial_{\theta} +
\frac{1}{2}
\partial_{hh} $ in $(h, \theta)$ coordinates, which starts at the
origin (corresponding to the point $x$). Let $\tau$ be the time of
the first exit from the following domain: $ \mathcal{D}^0 = \{
\theta < \theta (A_j) \} \cup \{ \theta \geq \theta(A_j); h > 0
\}$. Then define
\begin{equation} \label{p0}
p_0(x, dy) = \sum_{k: x,y \in \partial U_k} \sum_{ m \geq 0} {\rm
Prob}_x\{\theta(X_\tau) \in [\theta(y) + m  \int_{\partial U_k}
|\nabla H| dl , \theta(y+dy) + m  \int_{\partial U_k} |\nabla H|
dl ] \}.
\end{equation}
The summation over $k$ is needed  to account for the fact that $x$
and $y$ may both belong to the same edge $\gamma(A_i, A_j)$, in
which case they both belong to the boundaries of two cells, and we
need to consider two sets of $(h, \theta)$ coordinates.

The function $p_0(x,dy)$ is a stochastic transition function on
${\cal L}_p$,  and it can be considered as a stochastic transition
function on ${\cal L}^0\times S$. It clearly satisfies conditions
(A) and (B) preceding Lemma \ref{l1}.

For $x, y \in {\cal L}_p$, and for the stopping times
$\tau_n^{\varepsilon}$ defined in the previous section, let
$p_{\varepsilon}(x, dy)$ be the transition function for the chain
$X^{\varepsilon}_{\tau^{\varepsilon}_n}$ considered on  ${\cal
L}_p$:
\[
p_{\varepsilon}(x, dy) = {\rm Prob}_x \{
X^{\varepsilon}_{\tau_1^{\varepsilon}} \in [y , y+dy] \}~.
\]
Note that this definition is similar to that of $p_0(x, dy)$.
Here, however, we do not use the $(h, \theta)$ coordinates since
with small probability the process $X^{\varepsilon}_t$ may travel
outside of the domains $U_k$ for which $x \in \gamma(A_i,A_j)
\subseteq \partial U_k$ before time $\tau^{\varepsilon}_1$ (due to
the presence of the small diffusion term $X^{\varepsilon}_t$ may
go 'past' the saddle point $A_i$, thus traveling to one of the
neighboring domains before time $\tau^{\varepsilon}_1$).

\begin{lemma}
\label{l2} For any closed interval $I \subset \gamma(A_i, A_j)$,
and any bounded continuous function $f$ on ${\cal L}_p$,
\[
\lim_{\varepsilon \rightarrow 0} \int f(y) p_{\varepsilon}(x,dy) =
\int f(y) p_0(x,dy)~~~~{\rm uniformly~in ~}x \in I~.
\]
\end{lemma}
Notice that Lemma \ref{l2} implies the condition (C) for the chain
$Y^{\varepsilon}_n$. Before we start the proof of  Lemma \ref{l2}
we state and prove the
 following preliminary lemma.
\begin{lemma}
\label{X} Let $X^1_t$ and $X^2_t$ be the following two diffusion
processes on $ \R^d$ with infinitely smooth coefficients:
\[
dX^1_t = v(X^1_t)dt + a(X^1_t)dW_t + \varepsilon^2 v_1(X^1_t) dt +
\varepsilon a_1(X^1_t) dW_t~,
\]
\[
dX^2_t = v(X^2_t)dt + a(X^2_t)dW_t + \varepsilon^2 v_2(X^2_t) dt +
\varepsilon a_2(X^2_t) dW_t~,
\]
with $X^1_0 = X^2_0$. Suppose that for a certain constant $L$ the
following bound on the coefficients holds:
\[
|\nabla v^i| , |\nabla a^{ij}|, |v^i_1| , |v^i_2|, |a^{ij}_1|,
|a^{ij}_2| \leq L~,~~i,j = 1,...,d~,
\]
where $i$ and $j$ stand for the vector (matrix) entries of the
coefficients. Let $\lambda$ be the initial distribution for the
processes $X^1_0$ and $X^2_0$. Then for some constant $K = K(L)$
and for any $t, \eta
> 0$ we have
\[
{\rm Prob}_{\lambda} \{ \sup_{0 \leq s \leq t} |X_s^1 - X_s^2|
\geq \eta \} \leq \frac{(e^{Kt} - 1)\varepsilon^2}{\eta^2}~.
\]
\end{lemma}
\noindent  {\bf Proof:} Let us assume that $d = 1$ in order  to
avoid vector and matrix indices.  By Ito's formula, for any
stopping time $\tau \leq t$,
\[
\E_{\lambda}|X^1_{\tau} - X^2_{\tau}|^2 = \E_{\lambda}
\int_0^{\tau} 2(X^1_s - X^2_s)(v(X^1_s) - v(X^2_s) +
\varepsilon^2[v_1(X^1_s) - v_2(X^2_s)])ds +
\]
\begin{equation} \label{eql1}
\E_{\lambda} \int_0^{\tau} (a(X^1_s) - a(X^2_s) +
\varepsilon[a_1(X^1_s) - a_2(X^2_s)])^2 ds~.
\end{equation}
From the estimates on the coefficients and their derivatives it
follows that the expression in (\ref{eql1}) can be estimated as
follows
\begin{equation} \label{eql2}
\E_{\lambda}|X^1_{\tau} - X^2_{\tau}|^2 \leq
 K(L)( \int_0^{\tau}\E_{\lambda}|X^1_s -
X^2_s|^2 ds + \varepsilon^2 t).
\end{equation}
In particular, for $\tau = t$ we have \[  \E_{\lambda}|X^1_t -
X^2_t|^2 \leq
 K(L)( \int_0^t\E_{\lambda}|X^1_s -
X^2_s|^2 ds + \varepsilon^2 t).
\]
Let $R(t) = \E_{\lambda}|X^1_t - X^2_t|^2 + \varepsilon^2$. Then
$R(t) \leq K(L)\int_0^t R(s)ds,~R(0) = \varepsilon^2$. By
Gronwall's Lemma applied to R(t) we have
\[
\E_{\lambda}|X^1_t - X^2_t|^2 \leq \varepsilon^2(e^{K(L) t} -1)~.
\]
Define the stopping time $\tau = \min \{s: |X^1_s - X^2_s| \geq
\eta \} \wedge t$. Then, by (\ref{eql2})
\[
\eta^2 {\rm Prob}_{\lambda} \{ \max_{0 \leq s \leq t} ||X^1_s -
X^2_s| \geq \eta \} \leq \E_{\lambda}|X^1_{\tau} - X^2_{\tau}|^2
\leq
\]
\[
K(L)(\E_{\lambda} \int_0^{\tau}|X^1_s - X^2_s|^2 ds +
\varepsilon^2 t) \leq \varepsilon^2(e^{K(L) t} -1)~,
\]
which yields the lemma.
 \qed
 \\
 \\ {\bf Proof of Lemma \ref{l2}:}
Since the kernel $p_0(x, dy)$ is smooth in both variables and
satisfies condition (A) preceding Lemma \ref{l1}, to prove the
uniform weak convergence stated in Lemma \ref{l2} it is sufficient
to demonstrate that for an arbitrary closed interval $J \subset
\gamma(A_j, A_l)$ and an arbitrary $\delta >0$ there is
$\varepsilon_0 > 0$ such that
\begin{equation}
\label{est1} p_{\varepsilon}(x, J) > p_0(x, J) - \delta~~{\rm
for}~~{\rm all}~~ x \in I, ~~\varepsilon < \varepsilon_0~.
\end{equation}
Suppose that $I, J \subset \partial U_k$. (If $I$ and $J$ don't
belong to the boundary of the same cell, then  $p_0(x, J)$ is
equal to zero.) For the sake of simplicity of notation let us
assume that $I \subset \gamma(A_i, A_j)$ and $J \subset
\gamma(A_j, A_l)$, that is $I$ and $J$ belong to the adjacent
edges of ${\cal L}_p$. Without loss of generality we can assume
that $h
> 0$ in $U_k \cap V^{\varepsilon}$ for sufficiently small
$\varepsilon$. We can consider the process $X_t^{\varepsilon}$ in
$(h, \theta)$ coordinates in the following domain (see the picture
below)
\[
\mathcal{D}^{\varepsilon}= \{ \theta > \theta(A_i) \} \bigcap
\{|h| < \varepsilon^{\alpha_1 - \frac{1}{2}}\} \bigcap \left( \{
\theta < \theta(A_j) \} \bigcup \{ \theta \geq \theta(A_j), h>0
\}\right)~.
\]

\begin{figure}[htbp]
 \label{pic3}
  \begin{center}
    \begin{psfrags}
     \includegraphics[height=3in, width= 5.2in,angle=0]{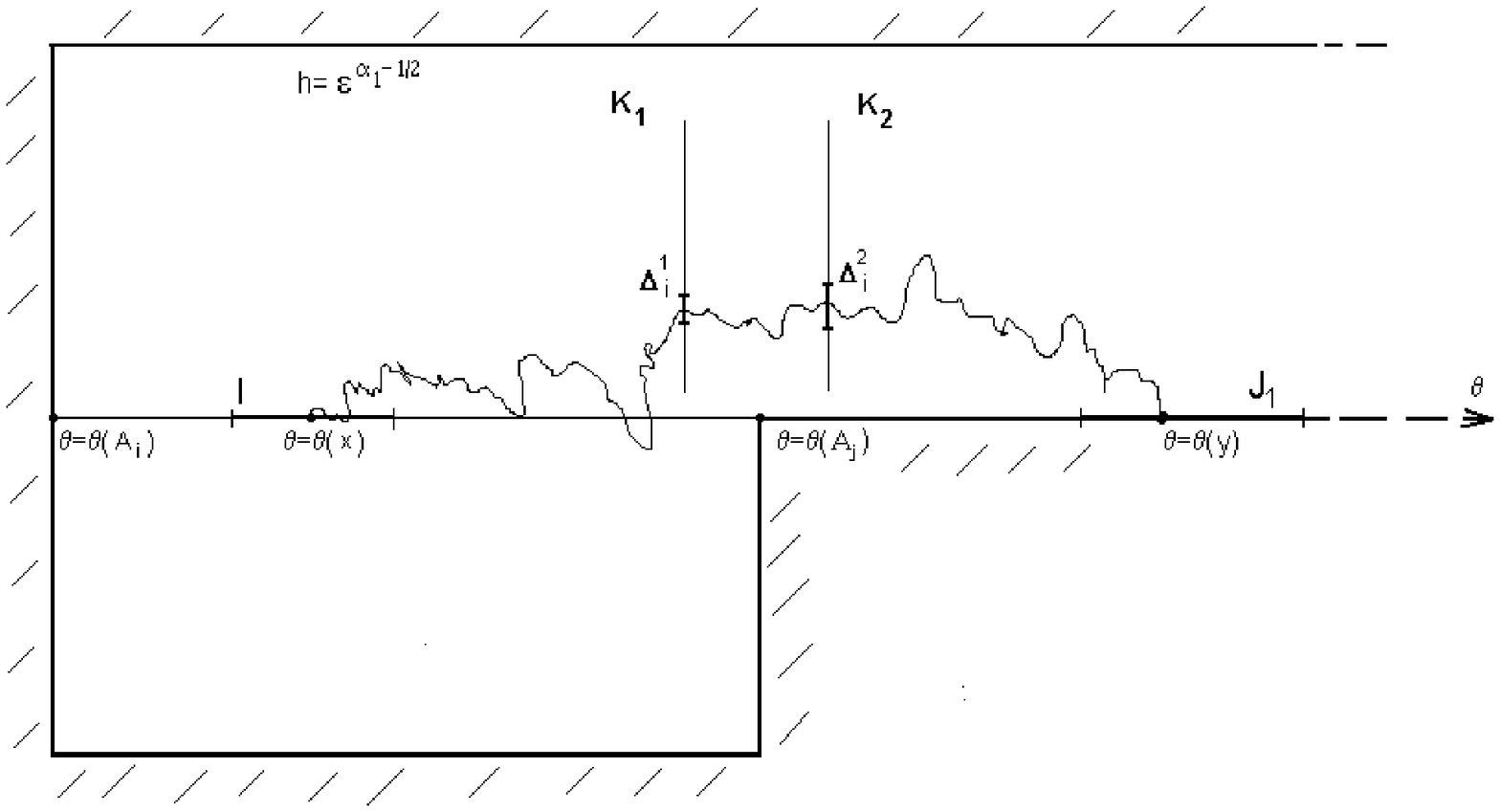}
    \end{psfrags}
  \end{center}
\end{figure}
As above we consider the process $X_t$ in $(h, \theta)$
coordinates in the domain \[
  \mathcal{D}^0 = \{
\theta < \theta (A_j) \} \bigcup \{ \theta \geq \theta(A_j); h > 0
\}~.
\]
Note that $p_{\varepsilon}(x, J)$ is estimated from below by the
probability that $X_t^{\varepsilon}$ leaves
$\mathcal{D}^{\varepsilon}$ through any of the copies of $J$
(which corresponds to $X_t^{\varepsilon}$ making a finite number
of rotations inside $U_k \cap V^{\varepsilon}$, and then leaving
$U_k$ through the segment $J$). Let $J_0, J_1, ...$ be the copies
of $J$ in $(h, \theta)$ coordinates ($J_{m+1} $ can be obtained
from $J_m$ by a shift by $\int_{\partial U_k} |\nabla H| dl$ along
the $\theta$ axis). For an initial point $x \in
\mathcal{D}^{\varepsilon}$ let $ \widetilde{p}_{\varepsilon}(x,
J_m)$ be the probability that $X_t^{\varepsilon}$ leaves the
domain $\mathcal{D}^{\varepsilon}$ through $J_m$. As stated above,
for $x \in I$ we have
\begin{equation}
\label{ineq1} p_{\varepsilon}(x,J) \geq \sum_{m=0}^{\infty}
\widetilde{p}_{\varepsilon}(x, J_m)~.
\end{equation}
Similarly, for $x \in \mathcal{D}^0$ let $ \widetilde{p}_0(x,
J_m)$ be the probability that $X_t$ leaves the domain
$\mathcal{D}^0$ through $J_m$. If $x \in I$ and $J$ belong to
different edges of $\partial U_k$ (as we have assumed) then by
(\ref{p0}) we have
\begin{equation}
\label{ineq2} p_0(x,J) = \sum_{m=0}^{\infty} \widetilde{p}_0(x,
J_m)~.
\end{equation}
From the definition of  $\widetilde{p}_0(x, J_m)$ it is easily
seen that the sum converges uniformly in $x$ for $x \in I$. In
order to prove (\ref{est1}) it is sufficient to demonstrate that
for each $m$ and each $\delta
> 0$ there is $\varepsilon_0 > 0$ such that
\begin{equation}
\label{est2}
 \widetilde{p}_{\varepsilon}(x, J_m) >
\widetilde{p}_0(x, J_m) - \delta~~{\rm for}~~{\rm all}~~ x \in I,
~~\varepsilon < \varepsilon_0~.
\end{equation}
Indeed,  for every positive $\delta$ we can select $m_0$ such that
$ p_0(x,J) < \sum_{m=0}^{m_0} \widetilde{p}_0(x, J_m) +
\frac{\delta}{2}$. By (\ref{est2}) we can  ensure that if we take
$\varepsilon_0$ to be sufficiently small then
\[
\widetilde{p}_0(x, J_m) < \widetilde{p}_{\varepsilon}(x, J_m) +
\frac{\delta}{2(m_0+1)}~~{\rm for}~~{\rm all}~~x \in I,
\varepsilon < \varepsilon_0, m \leq m_0 ~, \]
 which implies
(\ref{est1}) due to (\ref{ineq1}).

For the sake of simplicity of notation we shall only prove
(\ref{est2}) for $m =0$.

The generator of the process $X_t^{\varepsilon}$ is
\[
L^{\varepsilon}f = \frac{\varepsilon}{2} \Delta f + v \nabla f~,
\]
which in $(h , \theta)$ coordinates becomes
\[
L^{\varepsilon}f = \frac{1}{2}(f''_{hh}|\nabla H|^2 + \varepsilon
f''_{\theta \theta} |\nabla \theta|^2 + \sqrt{\varepsilon} f'_h
\Delta H + \varepsilon f'_{\theta} \Delta \theta) + f'_{\theta}
|\nabla H| |\nabla \theta|~.
\]
Dividing all of the coefficients of the generator by the same
function $|\nabla H| |\nabla \theta|$ amounts to a random time
change for the process  $X_t^{\varepsilon}$, which does not affect
any of the transition probabilities. We shall denote the
time-changed process with the generator
$\widetilde{L}^{\varepsilon}f = \frac{L^{\varepsilon}f}{|\nabla H|
|\nabla \theta|}$ also by $X_t^{\varepsilon}$. This process
 satisfies the equation \begin{equation}
\label{f1} dX_t^{\varepsilon} = (1,0) \sqrt{\frac{|\nabla
H|}{|\nabla \theta|}} dW^h_t + (0,1) \sqrt{\varepsilon}
\sqrt{\frac{|\nabla \theta|}{|\nabla H|}} dW^{\theta}_t +
(\frac{\sqrt{\varepsilon}}{2} \frac{\Delta H}{|\nabla
\theta||\nabla H|}, \frac{\varepsilon}{2} \frac{\Delta
\theta}{|\nabla \theta||\nabla H|} +1 ) dt~,
\end{equation}
while
\begin{equation} \label{f2}
dX_t =  (1,0)dW^h_t +(0, 1 ) dt~,
\end{equation}
where $W^h_t$ and $W^{\theta}_t$ are one dimensional Wiener
processes in $h$ and $\theta$ variables respectively. We can not
apply Lemma \ref{X} to (\ref{f1}) and (\ref{f2}) directly, as some
of the coefficients on the right hand side of (\ref{f1}) may be
unbounded near the saddle point $(h = 0, \theta = \theta(A_j))$.
To circumvent this problem we shall take a sequence of steps
(justified below), which will single out a small neighborhood of
the saddle point, where yet another coordinate system will be
considered.

Step 1. Let us take $\delta'> 0$ small enough, so that there exist
$0<h_1< h_2$, such that for any $\theta_1 \in [\theta(A_j) -
\delta', \theta(A_j)]$ the probability of the event that $X_t$
passes through the interval $K = \{h_1 \leq h \leq h_2, \theta =
\theta_1 \}$ before leaving the domain $ \mathcal{D}^0$ through
$J_0$ differs from $ \widetilde{p}_0(x, J_0)$ by less than
$\frac{\delta}{10}$ for any $x \in I$. Due to the smoothness of
the transition kernel of the process $X_t$, for some $\alpha < 1$
the interval $K$ can be replaced by any set $K'$ as long as $K'$
is contained in $K$ and has Lebesgue measure at least
$(h_2-h_1)\alpha$.

Step 2. If necessary make $\delta'$ from Step 1 smaller, so that
$|\widetilde{p}_0(A, J_0) - \widetilde{p}_0(B, J_0)| <
\frac{\delta}{10}$ whenever $|h(A) - h(B)| < \delta'$ and $A, B
\in [h_1, h_2] \times [\theta(A_j) - \delta', \theta(A_j) +
\delta']$.

Step 3. Take $\delta'' \leq \delta'$ and let $K_1 = \{ h_1 \leq h
\leq h_2, \theta = \theta(A_j) - \delta'' \}$, and $K_2 = \{ h_1
\leq h \leq h_2, \theta = \theta(A_j) + \delta'' \}$. Take
$\delta''$ sufficiently small so that whenever $A \in K_1$ the
process $X_t^{\varepsilon}$ starting at $A$ passes through the set
$\{ h(A) - \frac{(1-\alpha) \delta'}{2} < h < h(A) +
\frac{(1-\alpha) \delta'}{2},  \theta = \theta(A_j) + \delta'' \}$
before leaving $ \mathcal{D}^{\varepsilon}$ with probability at
least $1 - \frac{\delta}{10}$ for small enough $\varepsilon$.

Step 4. Let us split the set $K_2$ into intervals $\Delta^2_i~,
~i=1,...,r$ of length $\delta'$ (we can assume that each of the
intervals included the endpoints). Construct on $K_1$ the
intervals $\Delta^1_i~, ~i=1,...,r$ of length $\alpha \delta'$,
such that $h({\rm center}(\Delta^1_i)) = h({\rm
center}(\Delta^2_i))$. Let $\tau_{K_1}$ be the first time when a
process either exits $ \mathcal{D}^{\varepsilon}$($
\mathcal{D}^0$) or reaches $K_1$. Let us take $\varepsilon_0$ so
small that
\[
{\rm Prob}_x \{ X^{\varepsilon}_{\tau_{K_1}} \in \Delta_i^1 \}
\geq {\rm Prob}_x \{ X_{\tau_{K_1}} \in \Delta_i^1 \} -
\frac{\delta}{10r}~~{\rm for}~~{\rm all}~~ x \in I, ~~\varepsilon
< \varepsilon_0~.
\]

Step 5. Let us take $\varepsilon_0$ so small that
\[
 \widetilde{p}_{\varepsilon}(A, J_0) >
\widetilde{p}_0(A, J_0) - \frac{\delta}{10}~~{\rm for}~~{\rm
all}~~ A \in K_2, ~~\varepsilon < \varepsilon_0~.
\] \\
Assuming that the Steps 1 - 5 are valid let us prove (\ref{est2}).
By the Markov property
\[
\widetilde{p}_{\varepsilon}(x, J_0) \geq \sum_{i = 1}^r {\rm
Prob}_x \{ X^{\varepsilon}_{\tau_{K_1}} \in \Delta^1_i \} \min_{A
\in \Delta^1_i} \widetilde{p}_{\varepsilon}(A, J_0)~.
\]
By Steps 3, 5, and 2, the second factor on the right hand side can
be estimated as follows:
\[
 \min_{A \in
\Delta^1_i} \widetilde{p}_{\varepsilon}(A, J_0) \geq \min_{A \in
\Delta^2_i} \widetilde{p}_{\varepsilon}(A, J_0) -
\frac{\delta}{10} \geq
 \min_{A \in \Delta^2_i}
\widetilde{p}_0(A, J_0) - \frac{2\delta}{10} \geq
 \max_{A \in \Delta^1_i}
\widetilde{p}_0(A, J_0) - \frac{3\delta}{10} ~,
\]
while by Step 4
\[
{\rm Prob}_x \{ X^{\varepsilon}_{\tau_{K_1}} \in \Delta^1_i \}
\geq {\rm Prob}_x \{ X_{\tau_{K_1}} \in \Delta^1_i \} -
\frac{\delta}{10r}~.
\]
Combining the above inequalities and using Step 1 we obtain
\[
\widetilde{p}_{\varepsilon}(x, J_0) \geq \sum_{i = 1}^r({\rm
Prob}_x \{ X_{\tau_{K_1}} \in \Delta^1_i \} -
\frac{\delta}{10r})(\max_{A \in \Delta^1_i} \widetilde{p}_0(A,
J_0) - \frac{3\delta}{10}) \geq
\]
\[
 \sum_{i = 1}^r{\rm Prob}_x
\{ X_{\tau_{K_1}} \in \Delta^1_i \}\max_{A \in \Delta^1_i}
\widetilde{p}_0(A, J_0) - \frac{4\delta}{10} \geq \widetilde{p}(x,
J_0) - \frac{5\delta}{10}~,
\]
which implies (\ref{est2}).

It remains to justify the construction in Steps 1-5. The validity
of Steps 1 and 2 follows from the fact that the transition kernel
of the process $X_t$  is smooth. To justify Steps 4 and 5 it is
sufficient to consider both processes $X^{\varepsilon}_t$ and
$X_t$ in a compliment to a neighborhood of the saddle point, where
Lemma \ref{X} applies.

In order to justify Step 3 we note that by Morse Lemma in a
neighborhood $O_j$ of the saddle point $A_j$ there is a smooth
change of variables, such that in the new variables the stream
function is $H(x_1,x_2) = x_1x_2$, and the interior of $U_k$
corresponds to the first quadrant $x_1,x_2 >0$. In the new
variables the generator of the process $X^{\varepsilon}_t$, after
a random change of time, becomes $L^{\varepsilon}f =  \varepsilon
L_1f + v_1 \nabla f$, where $L_1$ is a differential operator with
first and second order terms, with bounded coefficients, and
$v_1(x_1, x_2) = (-x_1, x_2)$. We shall consider the operator
$L^{\varepsilon}$ in the domain
$\widetilde{\mathcal{D}}^{\varepsilon} = O_j \bigcap \{ x_1 > 0;
x_2> 0; x_1+x_2 > \varepsilon^{\frac{2}{3}} ; x_1x_2 <
\varepsilon^{\frac{1}{3}} \}$. Make a further change of variables
in $\widetilde{\mathcal{D}}^{\varepsilon}$:
\[
(x_1,x_2) \rightarrow (u,v) = (\frac{x_1x_2}{\sqrt{\varepsilon}},
x_2-x_1)~.
\]
In the new variables, after dividing all the coefficients of the
operator by $(x_1+x_2)$, which amounts to a random change of time
for the process, the operator can be written as
\begin{equation}
\label{star} L^{\varepsilon}f = M^{\varepsilon}f + \frac{\partial
f}{\partial v}~,
\end{equation}
where $ M^{\varepsilon}f$ is a differential operator with first
and second order terms. All the coefficients of $ M^{\varepsilon}$
can be made arbitrarily small in
$\widetilde{\mathcal{D}}^{\varepsilon}$ by selecting a
sufficiently small neighborhood $O_j$ of the point $A_j$, and then
taking $\varepsilon$ to be sufficiently small.

 The construction in Step 3 now follows from
Lemma \ref{X} by comparing the process whose generator is the
operator (\ref{star}) with the deterministic process with
generator $\frac{\partial f}{\partial v}$.
\\

\noindent {\bf Remark} To verify condition (D) (uniform mixing)
preceding Lemma \ref{l1} for the chain $Y^{\varepsilon}_n$ it is
sufficient to show (see \cite{Doob}, page 197) that there is an
integer $n \geq 1$, an interval $I \subset {\cal L}^0\times S$,
and a constant $c > 0$, such that
\begin{equation} \label{est101}
p^n_{\varepsilon}(x, dy) \geq c \lambda(dy), ~~{\rm and}~~  p^n_0
(x, dy) \geq c \lambda(dy), ~~{\rm for}~~ x \in {\cal L}^0 \times
S,~ y \in I~,
\end{equation}
where $\lambda(dy) $ is the Lebesgue measure on $I$. The proof of
estimate (\ref{est101}) is absolutely similar to that of Lemma
\ref{l2}.

\section{The Asymptotics of the Transition Times}
\label{se4}

In this section we shall study the asymptotics of the integral
$\int_{{\cal L}^0} \E_x \tau_1^{\varepsilon} d
\mu^\varepsilon(x)$, which enters in the expression (\ref{Dd}) for
the effective diffusivity.

We shall demonstrate the following:
\begin{equation} \label{fr1}
\E_x\tau^{\varepsilon}_1 \leq c \varepsilon^{-\frac{1}{2}}~~{\rm
for}~~{\rm all}~~x \in {\cal L}^0,
\end{equation}
\begin{equation} \label{fr2}
\lim_{\varepsilon \rightarrow 0} \varepsilon^{\frac{1}{2}} \E_x
\tau^{\varepsilon}_1 = f^0(x)~~{\rm uniformly}~~{\rm in}~~ x \in
I,
\end{equation}
where $f^0(x) \in C_b({\cal L}^0)$ is a positive function and $I$
is an arbitrary closed interval $I \subset \gamma(A_i, A_j)$. From
parts (a) and (b) of Lemma \ref{l1} it then follows that
\begin{equation} \label{aint}
\int_{{\cal L}^0} \E_x \tau_1^{\varepsilon} d \mu^\varepsilon(x) =
\varepsilon^{-\frac{1}{2}} ( \int_{{\cal L}^0} f^0(x) d \mu^0(x) +
o(1))~~{\rm as}~~\varepsilon \rightarrow 0~,
\end{equation}
where $\mu^0(x)$ is the invariant measure on ${{\cal L}^0}$ for
the kernel $p_0(x, dy)$, defined in Section \ref{se3}.

 The proof of formulas (\ref{fr1}) and (\ref{fr2}) will rely
on a sequence of lemmas stated below. We shall study separately
the probability of the event that the process ${X}^{\varepsilon}_t
$ starting form $x \in  \gamma(A_i, A_j)$ reaches $\partial
V^{\varepsilon}$ before time $\tau_1^{\varepsilon}$, and the
expectation of the time it takes for the process starting from
$\partial V^{\varepsilon}$ to reach ${\cal L}^0$.

Consider the process ${X}^{\varepsilon}_t $ together with the
process $X_t$, whose generator in $(h, \theta)$ coordinates is
$\frac{1}{2} \partial_{hh} + \partial_\theta$ in the domain
$\mathcal{D}^{\varepsilon}_1= \{ \theta(A_i) < \theta <
\theta(A_j); |h| < \varepsilon^{\alpha_1 - \frac{1}{2}} \} $.
 We follow the process $X_t$ till it exits $\mathcal{D}^{\varepsilon}_1$.
  Let $P_0(x, dh)$ be the corresponding transition kernel. Thus
$P_0(x, dh)$ coincides with a Gaussian distribution on
$-\varepsilon^{\alpha_1 - \frac{1}{2}} < h < \varepsilon^{\alpha_1
- \frac{1}{2}}$, and has two point masses at $h = \pm
\varepsilon^{\alpha_1 - \frac{1}{2}}$. Similarly let
$P_{\varepsilon}(x, dh)$ be the transition kernel for the process
${X}^{\varepsilon}_t $, which starts at $x \in \gamma(A_i, A_j)$
and is stopped at the time $\sigma_1^{\varepsilon}$. We have the
following:
\begin{lemma} \label{A}
For any continuous function $f: \R \rightarrow \R$, such that
$|f(h)| \leq 1+|h|$, there exists $c > 0$, such that
\begin{equation} \label{pt1}
\int |f(h)| P_{\varepsilon}(x, dh) < c~~{\rm for} ~~x \in
\gamma(A_i, A_j)~.
\end{equation}
Furthermore, for any closed interval $I \subset \gamma(A_i, A_j)$,
\begin{equation} \label{pt2}
\lim_{\varepsilon \rightarrow 0} \int f(h)( P_{\varepsilon}(x,
dh)- P_0(x, dh)) = 0~, ~~{\rm uniformly}~~{\rm in}~x \in I~.
\end{equation}
\end{lemma}
The proof of Lemma \ref{A} is completely similar to that of Lemma
\ref{l2}.

We introduce the following notation: $\tau^{V^{\varepsilon}}$ is
the first time the process ${X}^{\varepsilon}_t $ leaves
$V^{\varepsilon}$; similarly, $\tau^{U_k}$ and
$\tau^{V^{\varepsilon} \cap U_k}$ are the first instances when
${X}^{\varepsilon}_t $ leaves $U_k$ and ${V^{\varepsilon} \cap
U_k}$ respectively.

In order to estimate the probability that the process
${X}^{\varepsilon}_t $ starting from $x \in \gamma(A_i, A_j)$
reaches $\partial V^{\varepsilon}$ before time
$\tau^{\varepsilon}_1$ we shall need the following

\begin{lemma} \cite{FW} \label{B} There exists a constant $c > 0$, such
that
\[
\E_x \tau^{V^{\varepsilon}} \leq c \varepsilon^{2 \alpha_1-1} |\ln
\varepsilon|~~{\rm for}~~{\rm any} ~~x \in V^{\varepsilon}.
\]
\end{lemma}
This Lemma is the same as Lemma 4.7 of \cite{FW} (it must be
observed that the proof of Lemma 4.7 of \cite{FW} goes through for
any $A_{24} < \frac{1}{2}$).

In the event that ${X}^{\varepsilon}_{\sigma^{\varepsilon}_1} \in
V^{\varepsilon} \cap U_k$, after the stopping time
$\sigma_1^{\varepsilon}$ the process ${X}^{\varepsilon}_t $ may
exit $V^{\varepsilon} \cap U_k$ either through $\partial
V^{\varepsilon}$ or through ${\cal L}^0$. The next lemma estimates
the probability that the process exits the domain through
$\partial V^{\varepsilon}$.
\begin{lemma}
\label{C} There exists $c > 0$, such that for any $x \in
V^{\varepsilon} \cap U_k$ \begin{equation} \label{lce}  |{\rm
Prob}_x \{ {X}^{\varepsilon}_{\tau^{V^{\varepsilon}\cap U_k}}
\in
\partial V^{\varepsilon} \} - h(x) \varepsilon^{\frac{1}{2} -
\alpha_1}| \leq c \varepsilon^{\alpha_1} |\ln \varepsilon|~.
\end{equation}
\end{lemma}
\noindent {\bf Proof:} Let ${L}^{\varepsilon}$ be the generator of
the process ${X}^{\varepsilon}_t $ in the domain $V^{\varepsilon}
\cap U_k$. Then the probability in (\ref{lce}) is equal to the
solution $u(x)$ of the equation ${L}^{\varepsilon}u =0 $ in
$V^{\varepsilon} \cap U_k$ with the boundary conditions $u|_{H =
\varepsilon^{\alpha_1}} = 1$, $u|_{H = 0} = 0$. Let $u_1(x) = u(x)
- \frac{H(x)}{\varepsilon^{\alpha_1}}$. Then $u_1$ is the solution
of the equation ${L}^{\varepsilon}u_1  = -
{L}^{\varepsilon}\frac{H(x)}{\varepsilon^{\alpha_1}}$ with the
boundary conditions ${u_1}|_{\partial(V^{\varepsilon} \cap U_k )}
= 0$. By Lemma \ref{B}, since $\varepsilon^{-1}
{L}^{\varepsilon}H(x)$ is bounded uniformly in $\varepsilon$, the
solution $u_1$ is estimated as follows:
\[
|u_1| \leq c_0\varepsilon^{1-\alpha_1} \E_x \tau^{V^{\varepsilon}
\cap U_k} \leq c_1 \varepsilon^{\alpha_1}|\ln \varepsilon|~.
\]
This implies the statement of the lemma. \qed

Using the Markov property of the process ${X}^{\varepsilon}_t $
with respect to the stopping time $\sigma_1^{\varepsilon}$, we
obtain that  there is $c > 0$, such that for any $x \in {\cal
L}^0$ we have the following:
\begin{equation} \label{e1e}
 {\rm
Prob}_x \{ \tau^{V^{\varepsilon}} < \tau_1^{\varepsilon} \} \leq
\int_{-\infty}^{\infty} \,\, \sup_{\overline{x} \in
V^{\varepsilon}: h(\overline{x}) = \overline{h}} {\rm
Prob}_{\overline{x}}  \{
{X}^{\varepsilon}_{\tau^{V^{\varepsilon}\cap U(\overline{x})}}
\in
\partial V^{\varepsilon} \} P_{\varepsilon}(x, d \overline{h})
 \leq c
\varepsilon^{\frac{1}{2} - \alpha_1}~,
\end{equation}
where $U(\overline{x})$ is  the domain which contains
$\overline{x}$ (one of the domains $U_k$), and the second
inequality is due to Lemma \ref{C} and (\ref{pt1}).
 Furthermore, due to (\ref{pt2}) we can evaluate the asymptotics of the event $\{
\tau^{V^{\varepsilon}} < \tau_1^{\varepsilon};~
{X}^{\varepsilon}_{\tau^{V^{\varepsilon}}} \in \partial
V^{\varepsilon} \cap U_k \} $ as follows:
\begin{equation} \label{e2}
\lim_{\varepsilon \rightarrow 0} {\rm Prob}_x \{
\tau^{V^{\varepsilon}} < \tau_1^{\varepsilon};~
{X}^{\varepsilon}_{\tau^{V^{\varepsilon}}} \in \partial
V^{\varepsilon} \cap U_k \} \varepsilon^{\alpha_1 -\frac{1}{2}} =
\lim_{\varepsilon \rightarrow 0} \int_0^{\infty} \overline{h}
P_0(x, d\overline{h})~,~~{\rm uniformly}~~{\rm in}~~ x \in I.
\end{equation}
The next lemma allows us to estimate the expectation of the time
it takes for the process starting at $\partial V^{\varepsilon}$ to
return to ${\cal L}^0$.
\begin{lemma} \label{D}
For each of the domains $U_k$ there exists a constant $c_k > 0$,
such that
\begin{equation} \label{ed}
\lim_{\varepsilon \rightarrow 0} {\varepsilon^{1-\alpha_1}}{\E_x
\tau^{U_k}} = c_k~~{\rm uniformly}~~{\rm in}~~x \in \partial
V^{\varepsilon} \cap U_k.
\end{equation}
\end{lemma}
Note that Lemma \ref{D}, together with (\ref{e1e}) and (\ref{e2})
implies (\ref{fr1}) and (\ref{fr2}) since the expectation of the
time it takes for the process to reach $\partial V^{\varepsilon}$
can be estimated by Lemma \ref{B}. It remains to prove Lemma
\ref{D}.

We introduce notations and state several technical lemmas needed
for the proof of Lemma \ref{D}.

 Since $H(x) = 0$ on $\partial
U_k$,  we may assume without loss of generality that $H(x) > 0$
inside $U_k$ in a small neighborhood of $U_k$. Then there is a
region $V \subset U_k$, whose boundary consists of $\partial U_k$
and a level curve $\{ H(x) = H_0 \}$, and, by selecting a
sufficiently small $H_0$, we can ensure that each level set of $H$
in $V$ is connected and there are no critical points of $H$ in the
closure of $V$ other than on $\partial U_k$.

For $0 \leq H \leq H_0$ let us define the following functions:
\begin{equation} \label{cff}
a(H) = \int |\nabla H| dl, ~~b(H) = \int \frac{\Delta H}{|\nabla
H|} dl,~~ q(H) = \int \frac{1}{|\nabla H|} dl~,
\end{equation}
in each case the integration is over the level set $\{ H(x) = H, x
\in V \}$.  Let $r < \frac{H_0}{2}$ be a small constant, to be
specified later. Consider the function $f(H)$, which solves the
equation
\begin{equation} \label{eqf}
a(H) f''(H) + b(H) f'(H) = -q(H)~,
\end{equation}
with boundary conditions $f(0) = f(2r) =0$. While it not used here
explicitly, we note the fact that the operator in the left hand
side of (\ref{eqf}) after dividing it by the function $2\, q(H)$
becomes the generator of the limiting diffusion process on the
edge of the graph corresponding to the domain $U_k$ (cf \cite{FW}
and Section \ref{se5} of this article).

We need the following lemma, which will be proved in Section
\ref{se6}.
\begin{lemma}
\label{f} There is a function $g(r)$, which satisfies $\lim_{r
\rightarrow 0} g(r) = 0$, such that $|f'(H)| < g(r)$ for all $0 <
H < 2r$. Further, there is a constant $c > 0$ such that $|f''(H)|
< c|\ln H| $ and $|f'''(H)| < \frac{c}{H}$.
\end{lemma}

Let us select constants $\alpha_2$ and $\alpha_3$ such that
$\alpha_1 < \alpha_2 < \alpha_3 < \frac{1}{2}$. Define the subsets
$V^A$ and $V^B$ of $V$ as follows:
\[
V^A = \{x \in V; \varepsilon^{\alpha_2} < H(x) < r \},~~ V^B = \{x
\in V; \varepsilon^{\alpha_3} < H(x) < 2r \}~.
\]
Let $\tau^A$  be the first time  the process ${X}^{\varepsilon}_t
$ leaves $V^A$, similarly $\tau^B$ is the first time the process
leaves $V^B$.  Let $x_t$ be the deterministic process
\[
dx_t =  v(x_t) dt~,
\]
and let $T(x)$ be the time it takes the process $x_t$ starting at
$x$ to make one rotation along the level set, $T(x) = \inf_{t>0}
\{ x_t = x \}$. The next lemma shows that for times of order
$T(x)$ the process ${X}^{\varepsilon}_t $ is in a certain sense
close to the deterministic process $x_t$. The lemma is proved in
Section \ref{se6}.

\begin{lemma}
\label{close} For any $\delta > 0$ there is $\gamma > 0$ such that
\begin{equation} \label{tsx}
{\rm Prob}_x \{ \sup_{s \leq T(x)} |H({X}^{\varepsilon}_s) -
H(x_s)| > \varepsilon^{\frac{1}{2} -\delta} \} <
\varepsilon^{\gamma}~~{\rm for}~~{\rm all}~~x \in V^A.
\end{equation}
There exist $\delta' > 0$ and $\gamma > 0$ such that
\begin{equation} \label{tsn}
{\rm Prob}_x \{ \sup_{s \leq T(x)} |{X}^{\varepsilon}_s - x_s| >
\varepsilon^{\delta'} \} < \varepsilon^{\gamma}~~{\rm for}~~{\rm
all}~~x \in V^A.
\end{equation}
\end{lemma}

One of the main ingredients of the proof of Lemma \ref{D} is the
following lemma, which is  a particular case of the main result
(Theorem 2.3) of \cite{FW}.
\begin{lemma}\cite{FW}
\label{F} There is a constant $c_k > 0$, such that on each level
set $\{ H(x) = r, x \in V \}$ we have
\begin{equation}
\label{ef} \lim_{\varepsilon \rightarrow 0} \varepsilon \E_x
\tau^{U_k} = c_k (1+ g_1(r))~,
\end{equation}
the limit is uniform on each level set, and $g_1(r)$ satisfies
$\lim_{r \rightarrow 0} g_1(r) = 0$.
\end{lemma}
Lemma \ref{F} is different from Lemma \ref{D} in that the initial
point in (\ref{ef}) belongs to a fixed level set of $H$, while in
(\ref{ed}) the initial point is asymptotically close to $\partial
U_k$ as $\varepsilon \rightarrow 0$.
\\ \\ {\bf Proof of Lemma \ref{D}:} We shall demonstrate that
there exists a function $g(r)$, such that $\lim_{r \rightarrow 0}
g(r) = 0$, for which
\begin{equation}
\label{tt1} \E_x \tau^A \leq \varepsilon^{\alpha_1-1} g(r)~,
\end{equation}
uniformly in $x \in \{ H(x) = \varepsilon^{\alpha_1}, x \in V \}$.
Let us show that (\ref{tt1}) is sufficient to prove the lemma. As
in the proof of Lemma \ref{C}, from (\ref{tt1}) it follows that
\begin{equation}
\label{in01} |{\rm Prob}_x \{ H({X}^{\varepsilon}_{\tau^A}) = r \}
- \frac{\varepsilon^{\alpha_1}}{r} | \leq
c\frac{\varepsilon^{\alpha_1}}{r} g(r)~
\end{equation}
uniformly in $x \in \{ H(x) = \varepsilon^{\alpha_1}, x \in V \}$.
For $x \in \partial V^{\varepsilon} \cap U_k$ by the Markov
property
\begin{equation}
\label{equ} \E_x \tau^{U_k} = \E_x \tau^A +
\E_x(\E_{{X}^{\varepsilon}_{\tau^A}}\tau^{U_k};
H({X}^{\varepsilon}_{\tau^A}) =r) +
\E_x(\E_{{X}^{\varepsilon}_{\tau^A}}\tau^{U_k};
H({X}^{\varepsilon}_{\tau^A}) = \varepsilon^{\alpha_2})~.
\end{equation}
The first term on the right side of (\ref{equ}) is estimated from
above by $\varepsilon^{\alpha_1-1}g(r)$ due to (\ref{tt1}).  The
second term has the following asymptotics due to (\ref{in01}) and
Lemma \ref{F}
\[
|\E_x(\E_{{X}^{\varepsilon}_{\tau^A}}\tau^{U_k};
H({X}^{\varepsilon}_{\tau^A}) =r) - c_k \varepsilon^{\alpha_1-1}|
\leq {g}_2(r) \varepsilon^{\alpha_1-1}
\]
for sufficiently small $\varepsilon$, where ${g}_2(r)$ satisfies
$\lim_{r \rightarrow 0} {g}_2(r) = 0$. The last term on the right
side of (\ref{equ}) is estimated from above by $c
\varepsilon^{\alpha_2-1}$ due to Lemmas \ref{B} and \ref{C}, and
the repeated use of (\ref{equ}). Therefore
\[
|\E_x \tau^{U_k} - c_k \varepsilon^{\alpha_1-1} | \leq
\varepsilon^{\alpha_1 -1} g_3(r)
\]
for sufficiently small $\varepsilon$, and $g_3(r)$ which satisfies
$\lim_{r \rightarrow 0}g_3(r) = 0$. Since $r$ can be selected
arbitrarily small we obtain $\lim_{\varepsilon \rightarrow 0}
\frac{\E_x \tau^{U_k}}{\varepsilon^{\alpha_1-1}} = c_k$. It
remains to prove (\ref{tt1}).

Let $\sigma^B = \min \{ \tau^B, T(x) \}$. We shall prove that for
some $K > 0$ for all sufficiently small values of $r$
\begin{equation} \label{m2}
T(x) + \frac{ K \E_x
f(H({X}^{\varepsilon}_{\sigma^B}))}{\varepsilon} \leq \frac{ K
f(H(x))}{\varepsilon}
\end{equation}
for all $x \in V^A$. From (\ref{m2}) it follows that
\begin{equation} \label{es02}
\E_x \tau^A \leq \frac{K f(H(x))}{\varepsilon} ~~{\rm for}~~x \in
V_A.
\end{equation}
Due to the estimate on the derivative of $f$ from Lemma \ref{f}
for $x \in \{ H(x) = \varepsilon^{\alpha_1}, x \in V \}$ the right
side of (\ref{es02}) is estimated from above by $K g(r)
\varepsilon^{\alpha_1-1}$, which implies (\ref{tt1}). Now we need
to prove (\ref{m2}).

Applying Ito's formula to $f(H({X}^{\varepsilon}_t))$ we obtain
\[
\frac{1}{\varepsilon}( \E_x f(H({X}^{\varepsilon}_{\sigma^B})) -
f(H(x))) =\frac{1}{2} \E_x \int_0^{\sigma^B}(
f''(H({X}^{\varepsilon}_s))|\nabla H ({X}^{\varepsilon}_s)|^2 +
f'(H({X}^{\varepsilon}_s)) \Delta H ({X}^{\varepsilon}_s)) ds~,
\]
while from (\ref{eqf})
\[
T(x) = - \int_0^{T(x)} ( f''(H(x_s))|\nabla H (x_s)|^2 +
f'(H(x_s)) \Delta H (x_s)) ds~.
\]
Thus, what we want to show is that there is a constant $K_1$ such
that for all $x \in V^A$
\[
 \E_x|
\int_0^{\sigma^B} [ f''(H({X}^{\varepsilon}_s))|\nabla H
({X}^{\varepsilon}_s)|^2  -  f''(H(x_s))|\nabla H (x_s)|^2]ds | +
\]
\begin{equation} \label{bg}
 \E_x |\int_0^{\sigma^B}
[f'(H({X}^{\varepsilon}_s)) \Delta H ({X}^{\varepsilon}_s)) -
f'(H(x_s)) \Delta H (x_s))] ds| +
\end{equation}
\[
\E_x |\int_{\sigma^B}^{T(x)}  [ f''(H({X}^{\varepsilon}_s))|\nabla
H ({X}^{\varepsilon}_s)|^2 + f'(H({X}^{\varepsilon}_s)) \Delta H
({X}^{\varepsilon}_s))]ds |\leq K_1T(x)~.
\]
Since $f'$ and $\Delta H$ are bounded, and $K_1$ can be taken
arbitrarily large, it is sufficient to estimate only those of the
terms in (\ref{bg}) which contain the second derivative of $f$. By
Lemma \ref{close} we have ${\rm Prob}_x \{ \sigma^B < T(x) \} \leq
\varepsilon^{\gamma}$, and the second derivative of $f$ can be
estimated by Lemma \ref{f}. Therefore for the last term containing
$f''$ we have
\[
\E_x |\int_{\sigma^B}^{T(x)} f''(H({X}^{\varepsilon}_s))|\nabla H
({X}^{\varepsilon}_s)|^2 ds | \leq c T(x) |\ln
(\varepsilon^{\alpha_3})| {\rm Prob}_x \{ \sigma^B < T(x) \} \leq
cT(x)~.
\]
The estimate
\[
 \E_x|
\int_0^{\sigma^B} [ f''(H({X}^{\varepsilon}_s))|\nabla H
({X}^{\varepsilon}_s)|^2  -  f''(H(x_s))|\nabla H (x_s)|^2]ds |
\leq cT(x)
\]
follows from Lemma \ref{close} and the estimates on $f''$ and
$f'''$ of Lemma \ref{f}. This completes the proof of Lemma
\ref{D}. \qed
\\ \\ {\bf Proof of Theorem \ref{t1}:} The effective diffusivity
$D(\varepsilon)$ is related to the variance $d^{\varepsilon}$ of
the limit of the functional of the Markov chain, and to the
integral of the expectation of the transition times via formula
(\ref{Dd}). As shown in Lemma \ref{l2} and in the Remark following
it, Lemma \ref{l1} applies, and therefore there exists the limit
$d^0 = \lim_{\varepsilon \rightarrow 0} d^{\varepsilon}$. The
asymptotics of the integral $ \int_{{\cal L}^0} \E_x
\tau_1^{\varepsilon} d \mu^\varepsilon(x) $ is given by
(\ref{aint}). This completes the proof of the theorem. \qed

\section{The Case of the Open Channels}
\label{se5}

In this section it will be convenient to consider the process
$\widetilde{X}^{\varepsilon}_t$, which is the same as
$X^{\varepsilon}_t$, but only accelerated by the factor
$\frac{1}{\varepsilon}$, that is $\widetilde{X}^{\varepsilon}_t=
X^{\varepsilon}_{t/\varepsilon}$. This process satisfies the
equation
\[ d \widetilde{X}^{\varepsilon}_t = \frac{1}{\varepsilon} v(\widetilde{X}^{\varepsilon}_t)dt + dW_t,
\,\,\,\,\,\,\,\,\,\,\,\, \widetilde{X}^{\varepsilon}_t\in
\mathbb{R}^2. \] Note, that as a process on the torus,
$\widetilde{X}^{\varepsilon}_t$ is uniformly (in $\varepsilon$)
exponentially mixing. Following \cite{FW} we consider the finite
graph $G$ which corresponds to the structure of the level sets of
$H$ on the torus.

The graph $G$ is constructed as follows: we identify all the
points which belong to each connected component of each level set
of $H$. This way  each of the domains $U_k$, bounded by the
separatrices, gets mapped into an edge of the graph, while the
separatrices themselves get mapped into the vertices.  Let
$e(\widetilde{X}^{\varepsilon}_t)$ label the edge of the graph and
let  $H(\widetilde{X}^{\varepsilon}_t)$ be the coordinate on the
edge. Then the process $(e(\widetilde{X}^{\varepsilon}_t),
H(\widetilde{X}^{\varepsilon}_t))$ can be considered as a process
on the graph. It is proved in \cite{FW} (Theorem 2.2) that the
process $(e(\widetilde{X}^{\varepsilon}_t),
H(\widetilde{X}^{\varepsilon}_t))$ converges to a certain Markov
process on the graph with continuous trajectories, which is
exponentially mixing. We state the result here in less generality
than in \cite{FW}, but this is sufficient for our purposes.

\begin{theorem} \cite{FW} There is a Markov (diffusion) process $Y_t$ on the graph
$G$, which is exponentially mixing, and has continuous
trajectories, such that for any $T>0$, the process
$(e(\widetilde{X}^{\varepsilon}_t),
H(\widetilde{X}^{\varepsilon}_t))$ converges to $Y_t$ weakly in
$C([0,T],G)$.
\end{theorem}
We are now in the position to prove Theorem \ref{t2}.\\ {\bf Proof
of Theorem \ref{t2}}: The displacement of the process
$\widetilde{X}^{\varepsilon}_t$ in the direction $x_1$ (the
direction of the channels) is given by
\[ (\widetilde{X}^{\varepsilon}_t)^1=\frac{1}{\varepsilon} \int_0^t v_1(\widetilde{X}^{\varepsilon}_s)ds + W^1_t, \]
where $v_1$ is the $x_1$ component of the velocity field. Therefore,
\[ D^{11}(\varepsilon) =
\lim_{t \rightarrow \infty} \frac{ \E_{\lambda}
(\frac{1}{\varepsilon} \int_0^t
v_1(\widetilde{X}^{\varepsilon}_s)ds +
W^1_t)^2}{(\frac{t}{\varepsilon})}  =
 \frac{1}{\varepsilon}(\lim_{t\rightarrow\infty}
\frac{\E_{\lambda} (\int_0^t
v_1(\widetilde{X}^{\varepsilon}_s)ds)^2}{{{t}}}+o(1))=
\]
\[\frac{2}{\varepsilon}(\int_0^{\infty}\E_{\lambda}v_1(\widetilde{X}^{\varepsilon}_0) v_1(\widetilde{X}^{\varepsilon}_s)ds+ o(1))~, \]
where $\widetilde{X}^{\varepsilon}_0$ is distributed according to
the invariant (Lebesgue) measure $\lambda$ on $ \mathbb{T}^2$. For
a function $f\in C^{\infty}(\mathbb{T}^2)$, let $
\overline{f}(e,H)$, $(e,H) \in G$ be the function defined on the
graph, other than on the vertices, which is equal to the average
of $f$ over the corresponding connected component of the level set
of $H$
\[ \overline{f}(e,H)= \frac{\int_0^{T(x)}{f}(x_s)ds}{T(x)},\] where
$x_t$ is the solution of the deterministic equation
$dx_t=v(x_t)dt$, the initial point $x$ belongs to the level set,
and $T(x)$ is the time of one revolution around the level set. It
is easily seen that for any initial point $x$ which does not
belong to any of the separatrices of $H$ we have
\[\lim_{\varepsilon\rightarrow 0}\int_0^t \E_x f(\widetilde{X}^{\varepsilon}_s )ds=0, \,\,\,\,\,\,
\mbox{ if }\,\,\, \overline{f}(e, H)\equiv 0.\] Therefore,
\[\lim_{\varepsilon\rightarrow 0}\,\, \int_0^t
\E_{\lambda}
 v_1(\widetilde{X}^{\varepsilon}_0)v_1(\widetilde{X}^{\varepsilon}_s)\, ds
= \lim_{\varepsilon\rightarrow 0}\,\,  \int_0^t \E_{\lambda}
v_1(\widetilde{X}^{\varepsilon}_0)\overline{v}_1(e
(\widetilde{X}^{\varepsilon}_s), H(\widetilde{X}^{\varepsilon}_s))
\,\,  ds =\]
\begin{equation} \label{lo}     \int_0^t \E_{\mu} \overline{v}_1(Y_0)\overline{v}_1(Y_s)\, ds~,
\end{equation}
 where  $\mu$ is the measure on $G$, which is invariant
 for the process $Y_t$.

The integrals \[ \int_t^{\infty} \E_{\lambda}
v_1(\widetilde{X}^{\varepsilon}_0)v_1(\widetilde{X}^{\varepsilon}_s)\,
ds \] and \[   \int_t^{\infty} \E_{\mu}
\overline{v}_1(Y_0)\overline{v}_1(Y_s)\,ds \] can be made
arbitrarily small by selecting sufficiently large $t$ due to
uniform mixing of the processes $\widetilde{X}^{\varepsilon}_t $
and $Y_t$. Therefore,
\begin{equation} \label{ts}
\lim_{\varepsilon\rightarrow 0}\,\, \int_0^{\infty} \E_{\lambda}
v_1(\widetilde{X}^{\varepsilon}_0)v_1(\widetilde{X}^{\varepsilon}_s)\,
ds =   \int_0^{\infty}\E_{\mu}
\overline{v}_1(Y_0)\overline{v}_1(Y_s)\,ds,
 \end{equation}
  which
shows that the asymptotics for $D^{11}(\varepsilon)$ is as stated
in the theorem.

Now let us consider the asympotics for $D^{22}(\varepsilon)$. Note
that $ \overline{v}_2(e,H)\equiv 0$, thus the arguments leading to
(\ref{ts}) do not provide the asymptotics of
$D^{22}(\varepsilon)$. Let $P_1,\dots, P_n$ be those of the
separatrices of $H$ on the torus which, when unfolded onto the
plane, are non-compact. Let us select a point $A_i$ on each of
$P_i$. Let us introduce the sequence of stopping times $\tau_n$,
$n \geq 1$, which are the consecutive times when
$\widetilde{X}^{\varepsilon}_t $ makes the transition to a
different level set $P_i$. Thus
$\widetilde{X}^{\varepsilon}_{\tau_n} $ is a Markov chain on the
set $\{P_1, \dots ,P_n\}$. We can also consider the Markov chain
\[Z^{\varepsilon}_n=( \widetilde{X}^{\varepsilon}_{\tau_n}, \tau_n -\tau_{n-1}
, \Delta_n)\] on the extended phase space $ \{P_1, \dots P_n\}
\times \mathbb{R}_+ \times \mathbb{R}$. The third component
$\Delta_n$ is defined as follows: If
$\widetilde{X}^{\varepsilon}_{\tau_n} $ is considered on the
plane, then
\[\Delta_n= A^2(n) - A^2(n-1),\] where $A^2(n)$  is the $x_2$ coordinate of
the point corresponding to the separatrix containing the point
$\widetilde{X}^{\varepsilon}_{\tau_n} $.

Similarly we can introduce the stopping times $\eta_n$ for the process $Y_t$
on the graph, which are the consecutive times when $Y_t$ visits different  vertices
$Q_i=H(P_i)$ of $G$, corresponding to the unbounded separatrices of $H$.
Together with the Markov chain $Y_{\eta_n}$ we can consider the chain
\[ \widetilde{Z}_n=(Y_{\eta_n}, \eta_n-\eta_{n-1} , \widetilde{\Delta}_n )\] on
$ \{Q_1, \dots Q_n\} \times \mathbb{R}_+ \times \mathbb{R}$, where  $\widetilde{\Delta}_n $ is
defined the same way as $\Delta_n $.

Let $\mu ^{\varepsilon}$ be the invariant measure for the chain
$Z^{\varepsilon}_n$, and let $\widetilde{\mu}$ be the invariant
measure for the chain $\widetilde{Z}_n$. Let $f$ be the function
defined on the state space of the chain $Z^{\varepsilon}_n$, which
is equal to the third component: $f(x,\tau,\Delta)=\Delta$. The
function $ \widetilde{f}$ is defined the same way on $ \{Q_1,
\dots Q_n\} \times \mathbb{R}_+ \times \mathbb{R}$.

By the central limit theorem applied to the chain $Z_n^{\varepsilon}$, there is a number $d^{\varepsilon}$
such that
\[ \lim _{n\rightarrow\infty}\frac{\sum_{i=0}^n f(Z_i^{\varepsilon}) }{n}=N(0,d^{\varepsilon}).\] Similarly,
\[ \lim _{n\rightarrow\infty}\frac{\sum_{i=0}^n \widetilde{f}(\widetilde{Z}_i) }{n}=N(0,\widetilde{d}).\]

The effective diffusivity in the $x_2$ direction is then different form $d^\varepsilon$ by the factor
$\frac{1}{\varepsilon} \int \tau_1 d\mu^{\varepsilon}$,
 \[D^{22}(\varepsilon)=\frac{\varepsilon d(\varepsilon)}{\int \tau_1 d\mu^{\varepsilon}}. \]

From Theorem 2.2 in \cite{FW} and the uniform mixing of the Markov
chains $Z_n^{\varepsilon}$ and $\widetilde{Z}_n$ it easily follows
that \[d(\varepsilon)\rightarrow d\] and
\[\int \tau_1 d\mu^{\varepsilon}\rightarrow \int \eta_1 d\widetilde{\mu} . \]

This completes the proof of Theorem \ref{t2}. \qed

\section{Proof of the Technical Lemmas} \label{se6}
{\bf Proof of Lemma \ref{f}}: For any function $u \in C^{\infty}
(U_k) $ we have
\begin{equation}
\label{diff} \frac{d}{dH} \int u dl = \int \frac{ u \Delta
H}{|\nabla H |^2} dl + \int \frac{ \langle \nabla H, \nabla(
\frac{u}{|\nabla H |} ) \rangle }{|\nabla H| } dl~,
\end{equation}
where the integrals are over the level set $\{ H(x) = H, x \in V
\}$. In particular $b(H) = a'(H)$, and therefore equation
(\ref{eqf}) can be written as
\begin{equation}
\label{b1} (a(H)f'(H))' = -q(H)~.
\end{equation}
From the definition of the coefficients $a(H), b(H)$, and $q(H)$
it easily follows that
\begin{equation}
\label{abq} \lim_{H \rightarrow 0} a(H) = a_0 > 0~;~~ b(H) =
O(|\ln H|)~~{\rm as}~~H \rightarrow 0~;~~  q(H) = O(|\ln H|)~~{\rm
as}~~H \rightarrow 0~.
\end{equation}
Further, with the help of Morse Lemma and (\ref{diff}) it is
easily seen that
\begin{equation}
\label{bq}  b'(H) = O(\frac{1}{H})~~{\rm as}~~H \rightarrow 0~;~~
q'(H) = O(\frac{1}{H})~~{\rm as}~~H \rightarrow 0~.
\end{equation}
Let $H_m \in (0, 2r)$ be the point where $f(H)$ achieves its
maximum, thus $f'(H_m) = 0$. From (\ref{b1}) it follows that
\begin{equation}
\label{b2} f'(H) = \frac{ - \int_{H_m}^H q(s) ds}{a(H)}~.
\end{equation}
Thus, the estimate on the first derivative of $f$ stated in the
Lemma follows from (\ref{abq}). Rewrite (\ref{eqf}) as
\[
f''(H) = - \frac{q(H) +b(H) f'(H)}{a(H)}~.
\]
From (\ref{abq}) in now follows that $|f''(H)| \leq c |\ln H|$ for
some $c > 0$. Differentiating both sides of (\ref{eqf}) we obtain
\[
f'''(H) = - \frac{q'(H) + b'(H)f'(H) + b(H) f''(H) + a'(H) f''(H)
}{a(H)}~.
\]
The estimate on $f'''(H)$ now follows from the estimates on the
first two derivatives and from (\ref{abq}) and (\ref{bq}). This
completes the proof of the Lemma. \qed
\\
\\
{\bf Proof of Lemma \ref{close}}: The proof is based on the use of
Lemma \ref{X}. We can not however apply Lemma \ref{X} to the pair
of processes $X^{\varepsilon}_t$ and $x_t$ directly, since the
rotation time $T(x)$ grows logarithmically in $\varepsilon$ when
$x \in V^A$.

Let us establish the following property of the deterministic flow
$x_t$:

Let $0 =t_0 < t_1 < t_2<...< t_n $. Consider a process $y_t$,
which solves the equation \begin{equation} \label{flow} dy_t =
v(y_t)dt \end{equation}
 on each of the segments $[t_0, t_1), [t_1,
t_2),...,[t_{n-1}, t_n]$, with a finite number of jump
discontinuities $\lim_{t \rightarrow t_i +} y(t) - \lim_{t
\rightarrow t_i-} y(t) = p_i$, $ i = 1,...,n-1$. Then for any
positive $c$  there are positive $\kappa$ and $\delta'$ such that
under the conditions
\[
x_{t_0} = y_{t_0} \in V^A~; ~~ \sum_{i = 1}^{n-1} ||p_i|| <
\varepsilon^{\frac{1}{2} - \kappa} ~; ~~ t_n - t_0 \leq c| \ln
\varepsilon|
\]
we have
\begin{equation}
\label{form1a}
 \sup_{0 \leq t \leq t_n }||y_{t} - x_{t}|| < \varepsilon^{2 \delta'}~.
\end{equation}
Note that it is sufficient to establish the following: for any
pair of points $a_0, b_0$ such that $a_0 \in V^A$ and $||a_0 -
b_0|| \leq \varepsilon^{\frac{1}{2} - \kappa}$ we have
\begin{equation}
\label{one1}
 \sup_{0 \leq t \leq c |\ln \varepsilon| }||a_{t} - b_{t}||
  < \varepsilon^{2 \delta'+\kappa - \frac{1}{2}}||a_0 - b_0|| ~,
\end{equation}
where $a_t $ and $b_t$ are the solutions for the deterministic
flow (\ref{flow}). Let us take \begin{equation} \label{gr} \delta'
= \kappa = \frac{1}{4}(\frac{1}{2} - \alpha_2)~. \end{equation}
The time it takes for the trajectory of (\ref{flow}) to make one
rotation along the level set $\{ H(x) = H, x \in V \} $ is equal
to $T(x) = q (H(x))$, where $q(H)$ was defined in (\ref{cff}) and
is a smooth function for sufficiently small positive $H$, which
satisfies
\begin{equation}
q(H) = O(|\ln H|)~,~ q'(H) = O(\frac{1}{H})~{\rm as }~ H
\rightarrow 0~. \label{asym}
\end{equation}
The number of full rotations of the trajectory starting from $a_0$
is equal to $[\frac{t}{T(a_0)}]$ and the time it takes to make
these rotations is equal to  $[\frac{t}{T(a_0)}]T(a_0)$. It takes
$[\frac{t}{T(a_0)}]T(b_0)$ to make the same number of rotations
for the trajectory starting at $b_0$.

Due to (\ref{asym}) the difference is estimated as follows;
\begin{equation}
\label{ec1} |[\frac{t}{T(a_0)}] T(a_0) - [\frac{t}{T(a_0)}]
T(b_0)| \leq {\rm const} |\ln \varepsilon | \varepsilon^{-
\alpha_2} ||a_0 -b_0||~.
\end{equation}
Here we used the facts that $t \leq c |\ln \varepsilon|$ and that
$H(a_0) \geq \varepsilon^{\alpha_2}$. Now consider the images of
$a_0$ and $b_0$ under the flow (\ref{flow}) for time $t \leq
T(a_0)$. Using the reduction of the flow to a linear system in a
neighborhood of each of the saddle points (Hartman-Grobman Theorem
\cite{Pe}), it is easy to show that
\[
\sup_{0 \leq T(a_0)} ||a_t - b_t|| \leq {\rm const} ||a_0 - b_0||
\varepsilon^{-\alpha_2}~.
\]
Combining this with (\ref{ec1}) and with the fact that the speed
of motion in (\ref{flow}) is bounded, we obtain (\ref{one1}) with
$\delta'$ and $\kappa$ defined in (\ref{gr}). This in turn implies
(\ref{form1a}) as noted above.

Note that $0 < 2 \delta' < \frac{1}{2} - \kappa$ and that by
making $\kappa$ smaller (if necessary) we can satisfy $0 < \kappa
< \delta$, where $\delta$ is the same as in (\ref{tsx}). Observe
that for some $c> 0$ we have
\begin{equation}
\label{length}
  T(x) < c |\ln
\varepsilon|~~ {\rm for}~~{\rm  all}~~ x \in V^A.
\end{equation}
 Select the
points $0 =t_0 < t_1 < t_2<...< t_n = T(x)$ in such a way that
$\frac{\kappa}{2K}|\ln \varepsilon| \leq |t_{i+1} - t_i | \leq
\frac{\kappa}{K}|\ln \varepsilon|$ for $ i = 0,...,n-1$, where $K$
is the constant from Lemma \ref{X} (applied to the pair of
processes $X^{\varepsilon}_t$ and $x_t$). By (\ref{length}) there
is the estimate  $n \leq \frac{2 c K}{\kappa}$. Let
$y^{\varepsilon}_t$ be the piecewise continuous process, which is
defined by the conditions: $y^{\varepsilon}_{t_i} =
X^{\varepsilon}_{t_i}$ and $dy^{\varepsilon}_t =
v(y^{\varepsilon}_t) dt$ on $[t_i, t_{i+1})$, $i = 0,...,n-1$. By
Lemma \ref{X}
\begin{equation}
\label{y1}
 {\rm Prob}_x \{ \sum_{i = 0}^{n-1} \sup_{t \in [t_i,
t_{i+1})} || X^{\varepsilon}_{t} - y^{\varepsilon}_{t}|| >
\varepsilon^{\frac{1}{2} - \kappa} \} \leq (\frac{2 c
K}{\kappa})^3 \varepsilon^{\kappa}~.
\end{equation}
Due to continuity of $X^{\varepsilon}_t$  formula (\ref{y1})
provides an estimate on the sum of the jumps of the process
$y^{\varepsilon}_t$.  From (\ref{form1a}) it now follows that
\begin{equation}
\label{y2}
 {\rm Prob}_x \{|| x_{t} - y^{\varepsilon}_{t}|| >
\varepsilon^{2\delta'}\} \leq (\frac{2 c K}{\kappa})^3
\varepsilon^{\kappa}~.
\end{equation}
This, together with (\ref{y1}) implies (\ref{tsn}) for any $\gamma
< \kappa$. Since $H(x_t)$ is constant and $H(y^{\varepsilon}_t)$
is piecewise constant, we have
\[
{\rm Prob}_x \{ \sup_{t \leq T(x) } |H(X^{\varepsilon}_t ) -
H(x_t)| > \varepsilon^{\frac{1}{2} - \delta} \} \leq {\rm Prob}_x
\{ \sum_{i = 0}^{n-1} \sup_{t \in [t_i, t_{i+1})} ||
X^{\varepsilon}_{t} - y^{\varepsilon}_{t}|| >
\frac{\varepsilon^{\frac{1}{2} - \delta}}{\sup ||\nabla H|| } \}~.
\]
This, together with (\ref{y1}) and the condition that $\kappa <
\delta$  implies (\ref{tsx}) for any $\gamma < \kappa$. This
completes the proof of the Lemma. \qed
\\
\\
{\bf \large Acknowledgements}: I am grateful to Prof. S. Molchanov
for introducing me to this problem and for many useful
discussions.

\end{document}